\RequirePackage{fix-cm}
\documentclass[twocolumn]{svjour3}          
\smartqed  
\usepackage{graphicx}
\usepackage[numbers]{natbib} 
\usepackage{enumitem}
\usepackage{amsmath}
\usepackage{hyperref}
\newcommand{\mat}[1]{\mathbf{#1}}
\renewcommand{\tens}[1]{\boldsymbol{\mathsf{#1}}}
\newcommand{\diag}{\mathrm{diag}}
\newcommand{\n}{\vec{n}}
\newcommand{\dxdy}[2]{\frac{\partial #1 }{\partial #2 }}

\newcommand{\liq}{o}
\newcommand{\vap}{v}
\newcommand{\wat}{w}

\newcommand{\res}{F}
\newcommand{\ncell}{N}
\newcommand{\nblock}{M}
\newcommand{\npart}{N_p}
\newcommand{\ncomp}{N_c}
\newcommand{\nsparse}{d}
\newcommand{\noverlap}{b}

\newcommand{\acc}{A}
\newcommand{\flux}{G}
\newcommand{\src}{Q}

\newcommand{\cell}{\Omega}
\newcommand{\block}{\Omega^c}

\usepackage{xstring}
\usepackage{xcolor}

\begin{document}

\title{Accelerating Multiscale Simulation of Complex Geomodels by Use of Dynamically Adapted Basis Functions
}


\author{{\O}ystein S. Klemetsdal \and
        Olav M{\o}yner           \and
        Knut-Andreas Lie
}


\institute{{\O}.S. Klemetsdal \at
            Norwegian University of Science and Technology\\
            Tel.: +47 98 43 86 39\\
            \email{oystein.klemetsdal@ntnu.no}
            \and
            O. M{\o}yner \at
            Norwegian University of Science and Technology/SINTEF Digital, Norway \\
            \email{olav.moyner@sintef.no}
            \and
            K.-A. Lie \at
            Norwegian University of Science and Technology/SINTEF Digital, Norway \\
            \email{knut-andreas.lie@sintef.no}
}

\date{Received: date / Accepted: date}

\maketitle

 \begin{abstract}
   A number of different multiscale methods have been developed as a robust alternative to upscaling and as a means for accelerated reservoir simulation of high-resolution geomodels. In their basic setup, multiscale methods use a restriction operator to construct a reduced system of flow equations on a coarser grid, and a prolongation operator to map pressure unknowns from the coarse grid back to the original simulation grid. The prolongation operator consists of basis functions computed numerically by solving localized flow problems. One can use the resulting multiscale solver both as a CPR-preconditioner in fully implicit simulators or as an efficient approximate iterative linear solver in a sequential setting. The latter approach has been successful implemented in a commercial simulator. Recently, we have shown that you can obtain significantly faster convergence if you instead of using a single pair of prolongation-restriction operators apply a sequence of such operators, where some of the operators adapt to faults, fractures, facies, or other geobodies. Herein, we present how you can accelerate the convergence even further, if you also include additional basis functions that capture local changes in the pressure.  
\end{abstract}

\section{Introduction}

A problem is said to have a multiscale character if it is governed by parameters or mechanisms that act across a wide range of spatial or temporal scales. In many applications there is a clear separation between local processes taking place on a microscale and the macroscale behavior of the whole system. Flow in porous media does not have a clear scale separation, which can be seen from the elliptic Poisson equations, $\nabla\cdot(\tens K\nabla p) = q$, modelling single-phase, incompressible flow. Here, $\tens K$ represents rock permeability and has a multiscale structure in the sense that spatial variations are characterized by a wide spectrum of length scales. Over the past decade, we have seen a large number of so-called multiscale methods \citep{Hou1997,EfendievHou:09} that try to systematically account for small-scale variations and incorporate them in elliptic Poisson solvers formulated on a coarser scale. This is done by a pair of restriction and prolongation operators that restrict properties defined on a fine grid onto a coarse partition of the grid and map quantities from the coarse partition onto the fine grid, respectively. The prolongation operator is formed from a set of basis functions, computed numerically by solving localized versions of the original flow problem. 

Multiscale methods were originally developed as a robust alternative to upscaling, which only represents subscale variations through a homogenized coefficient $\tens{K}^\ast$. Multiscale methods not only produce reduced models for computing flow solutions on a coarser grid, but the basis functions represent a systematic means to propagate this solution to define approximate solutions on fine or intermediate scales. More important, these methods offer a natural way to accelerate multiphase flow simulations by reusing basis functions from one step to the next. For this, it is particularly important that the approximate pressure is monotone and that fluxes are mass conservative. Many multiscale methods have been proposed for incompressible flow on Cartesian grids or grids consisting of simplices, but only a few methods can provide realistic simulations of real hydrocarbon assets; see \citep{Lie2017ms-review} for a more detailed discussion. To this end, a multiscale method must handle complex flow physics and the rough, unstructured grids that are usually encountered in simulation models. To simulate real assets, one also has to incorporate models of wells and surface facilities, with accompanying strategies for controlling the injection and production of fluids, and implement robust nonlinear solvers and time-stepping strategies.

Realistic flow physics is described by two main classes of models on the reservoir scale. Black-oil type models lump the chemical hydrocarbon species into two pseudo-components (a light gas component and a heavier oil component) that can mix and appear in a liquid oleic and/or a gaseous phase at reservoir conditions. Compositional simulators represent  individual chemical species (or a wider range of pseudo-components) of hydrocarbons (and solvents) and use an equation of state to compute densities and fluid compositions at equilibrium states. In many cases, one can split the model equations into a flow equation for pressure and fluxes and a set of transport equations for saturations and compositions. These two are then solved consecutively in separate steps. First, the flow equation is solved while keeping saturations/concentrations fixed, and then the transport equations are solved with fixed pressures and fluxes \citep{watts1986compositional,trangenstein1989mathematical}. Flow equations typically have a certain elliptic character, whereas transport equations tend to be strongly hyperbolic \citep{bell1986conservation}. A sequential approach makes it easy to use different discretizations and solvers for the two types of equations, and, in particular, offers a natural way to incorporate multiscale methods developed for (elliptic) flow equations. Sequential solutions can reduce computational costs significantly but also lead to inaccurate solutions that fail to resolve the correct coupling between pressure and the other variables. This can to a large degree be mitigated by adding outer iterations over the flow equation and transport equation steps \citep{jenny2006adaptive}. This approach is used in state-of-the-art multiscale solvers and has produced good results for challenging black-oil \citep{Moyner16:spej,Kozlova16:spej,Kozlova16:ecmor,Lie2017ms-review,Bao2017} and compositional cases \citep{MoynerTchelepi2018}. 
In some cases, the best way to resolve couplings between flow and transport is nevertheless to use a fully implicit discretization and solve for all variables simultaneously. The resulting system of discrete equations is usually ill-conditioned and requires special solution strategies. One can use the sequential splitting procedure as an initial guess for the fully implicit system. Couplings between pressure and other variables are chiefly local \citep{Lacroix2003} and one can reduce the fully implicit system significantly by use of certain indicators for the splitting error \citep{MM18:ecmor}. So-called constrained-pressure-residual (CPR) methods \citep{Wallis1983, Wallis1985} constitute another popular approach. In CPR, one computes an inexpensive estimate to a reduced elliptic 'pressure' equation and uses this as a preconditioner for the full linearized system in an iterative Newton--Raphson type solution procedure. In early 2013, our group conducted the first tests with the multiscale finite-volume method \citep{Jenny2003} as a solver for the CPR step in a fully implicit black-oil simulator. With our prototype MATLAB simulator, we observed approximately 50\% reduction in computational costs for the two-phase SPE~10 benchmark \citep{spe10} and approximately 30\% reduction for a refined version of the three-phase SPE~1 benchmark \citep{spe1}. This was much less than we hoped for, and results were never published. \citet{Cusini15:jcp} later reported results from a more in-depth study. Herein, we go back and revisit this idea, this time with a more robust multiscale solver, the multiscale restriction-smoothed basis (MsRSB) method \citep{Moyner2016}, which has been developed to handle the type of complex, unstructured grids encountered in contemporary simulation models.

Reservoir models generally adapt to surfaces that describe both the external and internal geology. This gives highly deformed and skew cell geometries with high aspect ratios and large variations in interface areas. Experience shows that the coarse partition used in multiscale methods should try to preserve the stratigraphic architecture and the stratigraphy as closely as possible and thus adapt to faults, major fractures, and internal stratigraphic layering. Petrophysical properties usually form a cell-wise representation of volumetric elements like depositional environments, flow units, channels, and lobes from an underlying geological model. Preserving these may be important to maximize the accuracy of the reduced, coarse-scale equations. Generating a single coarse partition that takes all these constraints into account is very challenging and seldom possible. However, this is fortunately not necessary, since multiscale methods like the MsRSB method are usually both robust and relatively accurate for a wide variety of coarse partitions. On the other hand, \citet{Lie2017} recently demonstrated how one can improve the convergence rate for the pressure problem significantly by using multiple restriction/prolongation operators, where each operator targets specific geological features or provides increased resolution near wells. In this work, we investigate whether we can apply the same ideas to the CPR pressure preconditioner to improve convergence rates for the full system. We also discuss how to honor strong couplings between pressure and the other unknowns by use of dynamically adapted basis functions that capture local pressure changes. We demonstrate the new ideas on examples ranging from simple conceptual models to a WAG scenario posed on a realistic shallow-marine reservoir model.

\section{Model equations}

For a fluid with an aqueous phase and two hydrocarbon phases made up of $\ncomp$ components, we state conservation of mass for component $\beta$ as follows
\begin{equation}
    \label{eq:component-conservation}
    \begin{aligned}
      &\frac{\partial}{\partial t}\left(\phi\left[\rho_\wat S_\wat X_{\wat \beta}  + \rho_\liq S_\liq X_{\liq \beta}  + \rho_\vap S_\vap X_{\vap \beta}\right]\right)\\
      & \quad + \nabla \cdot \left(\rho_\wat X_{\wat \beta} \vec v_\wat + \rho_\liq X_{\liq \beta} \vec v_\liq + \rho_\vap X_{\vap \beta} \vec v_\vap \right) = q_\beta.
  \end{aligned}
\end{equation}
Here, $\rho_\alpha$ and $S_\alpha$ are the mass density and saturation of phase $\alpha$, $X_{\alpha\beta}$ is the mass fraction of component $\beta$ in phase $\alpha$, whereas $q_\beta$ is the source/sink term of component $\beta$. The Darcy phase velocities $\vec v_\alpha$ read
\begin{align}
      \vec v_\alpha & = -\frac{k_{r\alpha}}{\mu_\alpha} \tens K \left(\nabla p_\alpha - \rho_\alpha g \nabla z\right),
\end{align}
where $\mu_\alpha$ and $p_\alpha$ denote viscosity and pressure of phase $\alpha$, $\phi$ and $\tens K$ the rock porosity and permeability, $k_{r\alpha}$ models the reduced permeability experienced by one phase in the presence of another, $g$ is the gravity acceleration, and $z$ is the vertical coordinate. In most cases, the water component only appears in the aqueous phase, so that $X_{\wat,\beta} = 1$ if $\beta$ is the water component, and thermodynamic equilibrium is determined by the fugacity balance between liquid and vapor for each of the hydrocarbon components:
\begin{equation*}
  f_\beta^\liq(p,T,X_{\liq 1}, \dots, X_{\liq \ncomp}) =  f_\beta^\vap(p,T,X_{\vap 1}, \dots, X_{\vap \ncomp}).
\end{equation*}
The system is closed by assuming that the phases occupy the entire pore space, and that the mass fractions sum to unity for each phase:
\begin{equation*}
  S_\wat + S_\liq + S_\vap = 1, \quad \sum_{\beta=1}^{\ncomp} X_{\alpha\beta} = 1, \quad \alpha = \wat, \liq, \vap.
\end{equation*}
Finally, the phase pressures are related to each other through capillary pressure equations on the form
\begin{align*}
  p_\liq = p_\wat + P_{c\liq\wat}(S_\wat, S_\liq), \quad p_\vap = p_\liq + P_{c\vap\liq}(S_\liq, S_\vap).
\end{align*}

The type of multiscale methods considered herein were originally developed to solve Poisson type equations arising in incompressible flow models. Let us therefore explicitly state the corresponding flow equation for the special case of a two-phase model with an aqueous and a liquid phase. Assuming that each phase only consists of a single component, we can sum the conservation equations \eqref{eq:component-conservation} and use the fact that the saturations sum to unity to derive the following elliptic equation for the fluid pressure,
\begin{equation}
  \label{eq:incomp-elliptic}
  \begin{aligned}
      &\nabla \cdot \left(\frac{1}{\mu_\liq}\left[k_{r\liq} + \frac{\mu_\liq}{\mu_\wat}k_{r\wat}\right] \tens K \nabla  p_\liq \right) = q \\
      & \, - \nabla \cdot \left(\frac{k_{r\liq}}{\mu_\liq} \tens K \nabla P_{c\liq\wat} + \frac{\rho_\liq}{\mu_\liq}\left[k_{r\liq} + \frac{\rho_\wat}{\rho_\liq}\frac{\mu_\liq}{\mu_\wat} k_{r\wat}\right] g\tens K \nabla z\right).
  \end{aligned}
\end{equation}
This equation is coupled to a hyperbolic (or parabolic) saturation equation through the phase mobilities $\lambda_\alpha = k_{r\alpha}/\mu_\alpha$ and the capillary pressure function $P_{c\liq\wat}$. In the special case of linear relative permeabilities and equal viscosities, the viscous term does not depend on saturation and reduces to the single-phase case $\vec v = \tfrac{1}{\mu_\liq}\tens K \nabla p_\liq$. For fixed relative permeabilities, the coupling depends on the viscosity ratio. A similar relation holds for the gravity term: the coupling vanishes for linear relative permeabilities and equal viscosities and densities, but depends in the general case on the viscosity and density ratios.

\section{Discrete equations}
We introduce a grid with cells $\{\Omega_i\}_{i =1}^\ncell$, use backward Euler for the temporal discretization, and integrate over each cell in space using the divergence theorem and the midpoint rule to obtain the following discrete form of \eqref{eq:component-conservation}:
\begin{equation}
  \label{eq:residual-equations}
  \res_\beta^i = \acc^i_\beta + \sum_{j \in \mathcal{N}(i)}\flux^{i,j}_\beta - \src_\beta^i = 0, \quad i = 1, \dots, \ncell,
\end{equation}
where
\begin{align}
    \label{eq:accumulation}
    \acc_\beta^i & = \sum_{\alpha = \wat, \liq, \vap} A_{\alpha, \beta}^i, & (\text{Accumulation})\\
    \flux_\beta^{i,j} & = \sum_{\alpha = \wat, \liq, \vap} \flux_{\alpha,\beta}^{i,j}, & (\text{Intercell fluxes}) \\
    \src_\beta^i &  = \frac{\Delta t}{|\Omega_i|}(q_\beta)_i^{n+1}. & (\text{Sources/sinks})
\end{align}
Accumulation and intercell flux for each phase read 
\begin{align*}
    \acc_{\alpha,\beta}^i &= \left(\phi\rho_\alpha S_\alpha X_{\alpha\beta}\right)_i^{n+1} - \left(\phi \rho_\alpha S_\alpha X_{\alpha\beta}\right)_i^{n}, \\
   \flux_{\alpha,\beta}^{i,j} & = \frac{\Delta t}{|\Omega_i|} |\Gamma_{ij}| \left(\rho_\alpha X_{\alpha\beta} \vec v_\alpha \cdot \n \right)^{n+1}_{ij}.
\end{align*}
Subscript $i$ refers to values in cell $\Omega_i$ with bulk volume $|\Omega_i|$, subscript $ij$ refers to interface values at interface $\Gamma_{ij}$ with area $|\Gamma_{ij}|$, and superscript $n$ refers to the time step. The set $\mathcal{N}(i)$ consists of indices to all cells sharing a common interface with cell $i$. As is standard in reservoir simulation, we use a two-point approximation for the flux terms and single-point upstream evaluation for the interface mobilities, see e.g., \citep{lie2017mrstbook} for details, although other choices are equally possible.

To write the discrete system in matrix form, we first use a Schur complement reduction to remove well equations and closure relations local to each cell. This gives a system of $\ncomp \times \ncell$ equations, with $\ncell$ cell pressures and $(\ncomp-1)\times \ncell$ additional non-pressure variables as unknowns:
\begin{align*}
    \vec \res = (&\vec \res_1, \dots, \vec \res_{\ncomp})^T  & \\
    = (&\res_1^1, \dots, \res_1^\ncell, \dots, \res_{\ncomp}^1, \dots, \res_{\ncomp}^\ncell)^T, & \ncomp \times \ncell\\
    \vec x = (&\vec x_1, \vec x_2, \dots, \vec x_{\ncomp}) \\
    = (&p_1, \dots, p_\ncell, \quad & \ncell \\
       &x_2^1, \dots, x_2^\ncell, \dots, x_{\ncomp}^1, \dots, x_{\ncomp}^\ncell), \quad & (\ncomp-1) \times \ncell
\end{align*}
which we can write more compactly as $\vec \res = (\vec \res_p, \vec \res_s)^T$ and $\vec x = \left( \vec x_p, \vec x_s \right)^T$.
We can now write the system \eqref{eq:residual-equations} as $\vec \res(\vec x) = \vec 0$ and use Taylor expansion around $\vec x$,
\begin{equation*}
  \vec \res (\vec x + \Delta \vec x) = \vec \res(\vec x) + \mat J \Delta \vec x + \mathcal{O}(|\Delta \vec x|^2),
\end{equation*}
where $\mat J = \partial \mat \res/\partial \vec x$ is the Jacobian matrix of $\vec \res$. Neglecting higher-order terms, we obtain the Newton--Raphson method
\begin{equation*}
  \vec x^{k+1} = \vec x^k + \Delta \vec x^k, \quad
  -\mat J(\vec x^k) \Delta \vec x^k = \vec \res(\vec x^k).
\end{equation*}

\section{Preconditioning: the constrained pressure residual (CPR) method}
We write the linearized Newton system in block matrix form as
\begin{equation}
  \label{eq:fim-system}
  -\begin{bmatrix}
      \mat J_{pp} & \mat J_{ps} \\
      \mat J_{sp} & \mat J_{ss}
    \end{bmatrix}
    \begin{bmatrix}
      \Delta \vec x_p \\ \Delta \vec x_s
    \end{bmatrix}
 = \begin{bmatrix}
    \vec \res_p \\
    \vec \res_s
  \end{bmatrix}.
\end{equation}
In most cases, \eqref{eq:fim-system} contains so many unknowns that we must use an iterative linear solver, whose convergence depends highly on the spectrum of the matrix. There are two main factors that contribute to impede convergence: the mixed elliptic--hyperbolic character of the flow equations, which makes pressure a strong variable, and large aspect ratios and variations in rock properties, which give high condition numbers. An effective preconditioning strategy \citep{Lacroix2003} is thus required to ensure the effectiveness of the iterative solver. Herein, we apply a so-called constrained-pressure-residual (CPR) method \citep{Wallis1983,Wallis1985,Gries2014}, which relies on an inexpensive estimate of the pressure update $\Delta \vec x_p$ as an initial guess for the solution to the entire system.


The first step in a CPR method is to decouple \eqref{eq:fim-system} into a system on the form
\begin{equation*}
  -\begin{bmatrix}
    \mat J^*_{pp} & \mat J^*_{ps} \\
    \mat J_{sp} & \mat J_{ss}
  \end{bmatrix}
  \begin{bmatrix}
    \Delta \vec x_p \\ \Delta \vec x_s
  \end{bmatrix}
  =
  \begin{bmatrix}
    \vec \res^*_p \\
    \vec \res_s
  \end{bmatrix},
\end{equation*}
where the matrix product $(\mat J_{pp}^*)^{-1}\mat J_{ps}^*$ is sufficiently small so that the solution to $\mat J_{pp}^*\Delta \vec x_p = \vec \res_p^*$ is a reasonable approximation to the pressure. In addition to reducing the coupling strength, the decoupling should ultimately \citep{Stuben2007}
\begin{itemize}
\item reduce the condition number of the overall system;
\item reduce the condition number of the block $\mat J_{pp}$;
\item render $\mat J_{pp}^*$ and $\mat J_{ss}^*$ as M-matrices;
\item be computationally inexpensive to perform.
\end{itemize}
Popular decoupling strategies include alternate-block factorization (ABF) \citep{Bank1989}, IMPES and quasi-IMPES (IMplicit Pressure, Explicit Saturation) \citep{Coats2000, Lacroix2000}, and dynamic row sum (DRS) \citep{Gries2014}. Herein, we follow the IMPES and quasi-IMPES approaches. Like most of the decoupling methods described in the literature, these two methods utilize the fact that couplings between pressure and the other quantities are chiefly local \citep{Lacroix2003}. Algebraically, this implies that the diagonal terms of the coupling block $\mat J_{ps}$ dominate the off-diagonal terms. Thus, since we only seek an approximation to the pressure update, it is sufficient to reduce diagonal coupling terms of $\mat J_{ps}$. The IMPES decoupling amounts to finding weights $\vec w_\beta = (w_{\beta,1}, \dots, w_{\beta,\ncell})^T$ so that the derivatives of the accumulation terms with respect to non-pressure variables cancel:
\begin{equation*}
  \sum_{\beta=1}^{\ncomp}  \mat W_\beta \dxdy{\mat \acc_\beta}{\vec x_k} = \mat 0, \quad k = 2, \dots, \ncomp.
\end{equation*}
Here, $\mat W_\beta$ is a diagonal matrix with $(\mat W_\beta)_{ii} = w_{\beta,i}$, and $\mat \acc_\beta = (\acc_\beta^1, \dots, \acc_\beta^\ncell)$ is the vector of accumulation terms for the $\beta$-th component, as defined in \eqref{eq:accumulation}. The quasi-IMPES approach aims to eliminate all couplings between pressure and non-pressure variables in the same cell:
\begin{equation*}
  \sum_{\beta=1}^{\ncomp}  \mat W_\beta\, \diag \left(\dxdy{\mat \res_\beta}{\vec x_k}\right) = \mat 0, \quad k = 2, \dots, \ncomp,
\end{equation*}
where $\diag(\partial \mat \res_\beta/\partial \vec x_k)$ refers to the diagonal of the Jacobian block $\partial \mat \res_\beta/\partial \vec x_k$. In either case, the decoupling procedure amounts to solving a system on the form
\begin{align*}
  &\begin{bmatrix}
    \mat M_{2,1}^T & \dots & \mat M_{2,\ncomp}^T \\
    \vdots & \ddots & \vdots \\
    \mat M_{\ncomp,1}^T & \dots & \mat M_{\ncomp,\ncomp}^T \\
    \mat I & \dots & \mat I
  \end{bmatrix}
  \begin{bmatrix}
    \vec w_1 \\ \vdots \\ \vdots \\ \vec w_{\ncomp}
  \end{bmatrix}
  =
  \begin{bmatrix}
    \mat 0 \\ \vdots \\ \mat 0 \\ \mat 1
  \end{bmatrix}, \\
  & \quad \text{where} \quad
  \mat M_{k,\beta} = \begin{cases}
    \dxdy{\mat \acc_\beta}{\vec x_k} & \text{IMPES}, \\
    \diag\left(\dxdy{\mat \res_\beta}{\vec x_k}\right) & \text{quasi-IMPES},
  \end{cases}
\end{align*}
and then premultiply $\mat J$ and $\mat \res$ by
\begin{equation*}
  \mat W =
  \begin{bmatrix}
    \mat W_{1} & & \dots &\mat W_{\ncomp} \\
    \mat 0 & \mat I & \dots & \mat 0 \\
    \vdots & & \ddots & \vdots \\
    \mat 0 &  & \dots  &\mat I
  \end{bmatrix},
  \quad \text{where} \quad (\mat W_\beta)_{ii} = w_{\beta,i}.
\end{equation*}
The CPR method presumes that the solution to the elliptic pressure equation dominates the solution to the full system. When this is not the case, CPR may not be beneficial at all. \citet{Gries2014} describes a heuristic approach to turn off two-stage preconditioning when the coupling is stronger than a given threshold.

\section{Multiscale methods}

In this section, we briefly describe the framework of algebraic multiscale methods, before we discuss how to apply such methods as pressure solvers in the CPR method. Consider a discrete and linearized pressure equation on the form
\begin{equation}
  \label{eq:pressure_system}
  \mat A \vec p = \vec q,
\end{equation}
defined over a fine grid $\{\cell_i\}_{i = 1}^\ncell$ that incorporates all details of the geological model. Multiscale methods start from a coarse partition $\{\block_i\}_{i = 1}^M$ of the fine grid consisting of $\nblock<\ncell$ blocks defined so that each cell $\cell_i$ in the fine grid belongs to a single block $\block_j$ in the coarse grid. We then associate basis functions that map the pressure degrees of freedom on the coarse grid to degrees of freedom on the fine grid. We can collect the basis functions in a prolongation operator $P :\{\block_j\} \rightarrow \{\cell_i\}$ and express it in matrix form as an $\ncell \times \nblock$ matrix $\mat P$. Element $i,j$ of $\mat P$ is the value of the $j$th basis function in the $i$th cell. Given a coarse-scale approximation $p_c$, we can now compute an approximation to the fine-scale pressure by prolongation of $p_c$ back to the fine grid, $\vec p \approx \mat P \vec p_c$. Analogously, we define a restriction operator $R: \{\cell_i\} \rightarrow \{\block_j\}$ that maps quantities from the fine grid onto the coarse blocks. We write this operator in matrix form as an $\nblock\times \ncell$ matrix $\mat R$. To form a reduced linear system on the coarse grid, we substitute $\vec p$ by $\mat P \vec p_c$ in \eqref{eq:pressure_system} and multiply by the restriction operator from the left on both sides of the equation,
\begin{equation}
  \label{eq:pressure_system_reduced}
  (\mat R \mat A \mat P)\vec p_c = \mat R \vec q, \quad \text{ or } \quad \mat A_c \vec p_c = \vec q_c.
\end{equation}
To be feasible, any pair of operators $(P,R)$ should fulfill the following three requirements \citep{Lie2017}:
\begin{enumerate}
\item Both operators are defined over the same non-overlapping coarse partition of the fine grid. Each column of $\mat P$ constitutes a basis function, and is associated with a coarse grid block.
\item The support of each basis function is compact and must contain the associated coarse block. In other words, the support region $S_j$ of the $j$th basis function must satisfy the following conditions: $\block_j \subset S_j \subset \cup_{k = 1}^\nblock \block_k$
\item The columns of $\mat P$ and the rows of $\mat R$ must form partitions of unity over the fine grid; that is, each row of $\mat P$ has unit row sum and each column of $\mat R$ has unit column sum.
\end{enumerate}
Specific multiscale methods differ in the way they define the prolongation and restriction operators. This type of method was first proposed by \citet{Hou1997} and has later been developed in many different directions, see \citep{Lie2017ms-review} for a review focused on methods that could be applicable to practical reservoir simulation. The main contenders used to be the multiscale finite-volume (MsFV) method \citep{Jenny2003} and the multiscale mixed finite-element (MsMFE) method \citep{Aarnes2006}, but these have later been superseded by the multiscale restriction-smoothed basis (MsRSB) method \citep{Moyner2016}. The mixed method is formulated for an extended mixed hybrid system that also contains unknowns for phase fluxes and pressures at cell interfaces.

\begin{figure*}[h]
    \centering
    \includegraphics[width=0.95\textwidth]{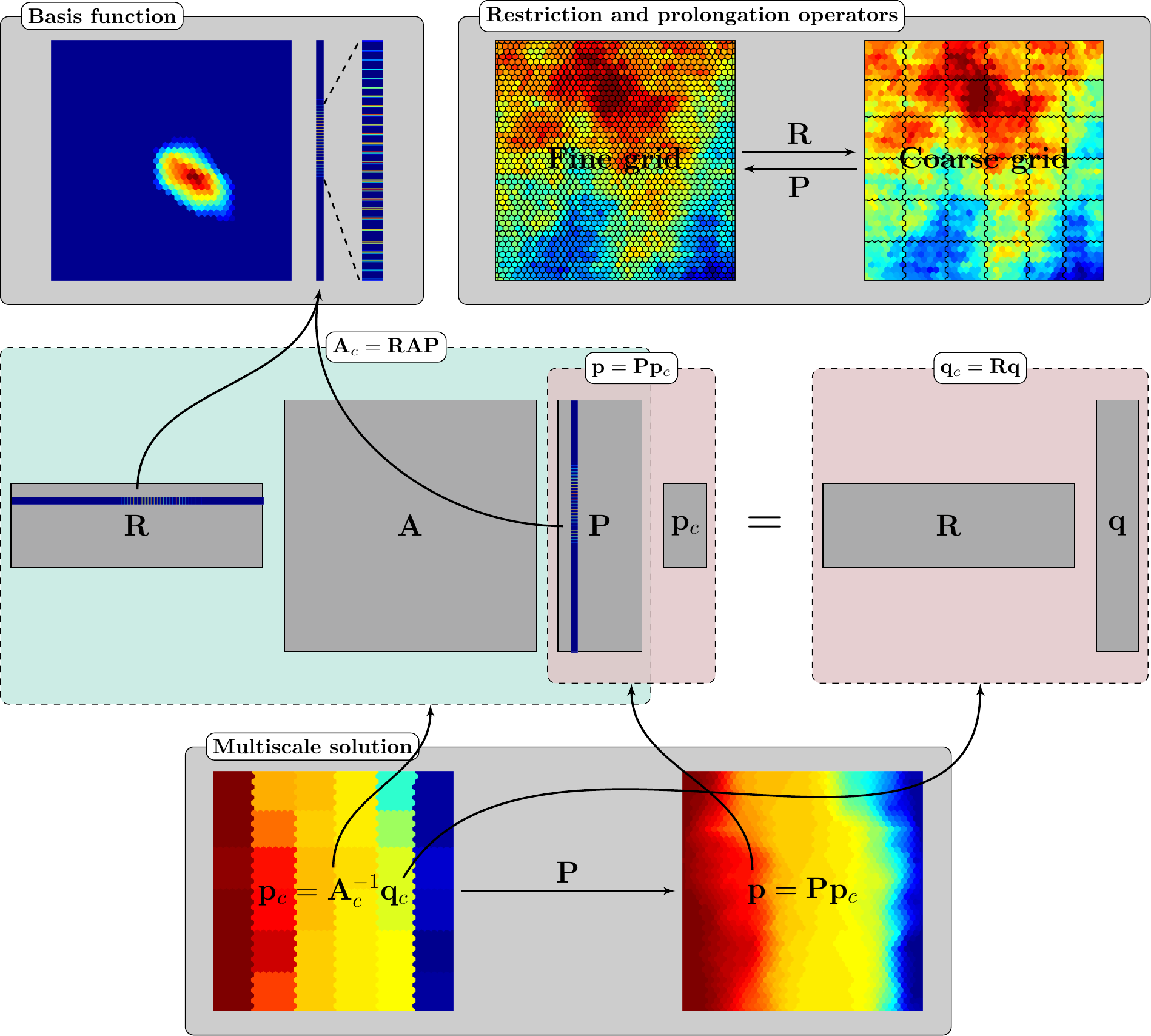}
    \caption{Illustration of the multiscale solution procedure. The middle section shows how the full system reduces to $\mat A_c \mat p_c = \mat q_c$. In the upper right corner, we see how the restriction and prolongation operators map between the two scales. In the upper left corner, one of the basis functions is depicted in physical space and in vector form as it appears in the restriction and prolongation operators as a row or column. The solution procedure is shown at the bottom, where we first obtain the coarse-scale pressure by solving the reduced system, and then prolongate it onto the fine grid to obtain the fine-scale pressure solution.}
\label{fig:multiscale}
\end{figure*}

How one should interpret \eqref{eq:pressure_system_reduced} depends on the choice of operators. For the restriction operator, a natural choice is to set $\mat R = \mat P^T$, which corresponds to a Galerkin coarse-scale discretization. This does not give a locally mass-conservative scheme on the coarse scale, and is consequently not a good choice in a standalone solver for the pressure equation. If fine-scale fluxes are needed, it is better to use a finite-volume operator that sums the entries of the corresponding cells for each block. Herein, we use multiscale as a preconditioner \citep{Wang14:jcp} and ensure mass conservation by converging the full system. Galerkin restriction is then the preferable choice, since it preserves the (hopefully) symmetric and M-matrix properties of the fine-scale pressure matrix. Figure~\ref{fig:multiscale} depicts the multiscale solution procedure, which we can interpret as a two-level algebraic multigrid method, albeit with a very large coarsening ratio.

\subsection{Iterative multiscale multibasis solver}
\label{sec:iterative-ms-mb}

As in standard multigrid theory \citep{trottenberg2000multigrid}, the multiscale solution typically resolves global low-frequency errors quite effectively, but contains local high-frequency errors due to the localization introduced to define basis functions. It is therefore natural to cast the multiscale method in an iterative framework \citep{msfvi,Wang14:jcp} and combine the multiscale solver with a smoothing step that aims to remove high-frequency errors. One iteration then consists of two steps:
\begin{align*}
  \label{eq:iterative_multiscale}
  \vec x^{k+1/2} & = \vec x^k + S(\vec A, \vec q - \mat A \vec x^k), \\
  \vec x^{k+1} & = \vec x^{k+1/2} + \mat P \mat A_c^{-1} \mat R (\vec q - \mat A \vec x^{k+1/2}),
\end{align*}
where $S(\mat A, \vec b)$ denotes a function that performs one or more smoothing operations, e.g., incomplete LU factorization. 

The initial convergence rate of such a two-stage multiscale preconditioner is typically satisfactory, but deteriorates when the error is made up of intermediate-frequency modes that are neither reduced by the smoother nor the multiscale step. When used as pressure solver in a CPR method, convergence to a strict tolerance is not necessary, and we are instead interested in reducing the error as much as possible using only one or a few iterations. \citet{Lie2017} previously showed that the initial reduction can be accelerated significantly if we instead of using basis functions associated with a single coarse partition use a sequence of prolongation/restriction operators associated with different coarse partitions that each targets specific features of the flow field and/or the geological model. In particular, we define a set of prolongation operators $\{P^1, \dots, P^{\npart}\}$ and restriction operators $\{R^1, \dots, R^{\npart}\}$ that each corresponds to a (unique) coarse grid $\{\block_i\}_{i = 1}^{\nblock_\ell}$ and satisfies requirements 1 to 3 above. This yields an iterative multiscale multibasis method on the form
\begin{equation}
    \label{eq:multiscale-multibasis}
    \begin{aligned}
        \vec x^{k+(2\ell -1)/2\npart} & = \vec x^{k+(\ell-1)/\npart} \\
        & \quad + S^\ell(\vec A, \vec q - \mat A \vec x^{k+(\ell-1)/\npart}), \\
        \vec x^{k+\ell/\npart} & = \vec x^{k+(2\ell -1)/2\npart} \\
        & \quad + \mat P^\ell (\mat A^\ell_c)^{-1} \mat R^\ell (\vec q - \mat A \vec x^{k+(2\ell -1)/2\npart}).
    \end{aligned}
\end{equation}
Employing a nested sequence of coarse partitions to capture different error modes is the standard approach in algebraic multigrid methods. In our case, the coarse partitions are all defined relative to the original fine grid and are not necessarily nested.

\begin{figure*}[h!]
    \centering
    \includegraphics[width=0.95\textwidth]{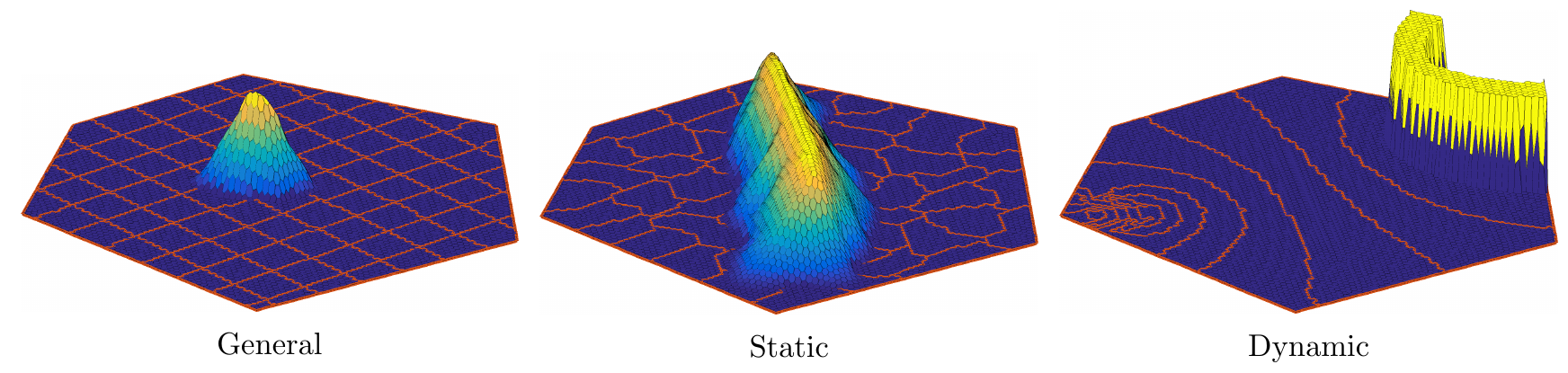}
    \caption{Examples of different types of basis functions. General basis functions are defined over a structured partition, the static partition traces a major fracture inside the domain, and the dynamic partition is based on $\Delta p$. The first two types are constructed using MsRSB, whereas the dynamic partition uses constant basis functions.}
    \label{fig:multiscale-bases}
\end{figure*}

Careful analysis of many simulation cases have lead us to the following three observations:
\begin{enumerate}[label=(O\arabic*)]
\item Multiscale operators defined for partitions that cover the entire domain evenly will usually resolve the global pressure (as given by the CPR pressure equation) quite accurately, but may introduce local errors that stem from localization of the basis functions. Such errors typically arise when the support of the basis function contains barriers or regions with distinctively different flow properties, or are the result of significant changes in the drive mechanisms.
\item Local errors introduced by geological variations, or in near-well regions, can be more effectively reduced by coarse partitions that adapt to these spatial features. 
\item Differences between the CPR pressure and the true pressure tend to arise because of inexact decoupling of the pressure equation. This effect increases with the coupling strength and typically gives temporal errors located near propagating displacement fronts. By adapting partitions to dynamic changes in the pressure update, or saturations and/or mass fractions, one can hope to obtain multiscale operators that resolve this coupling better.
\end{enumerate}
With these observations in mind, we herein consider three types of multiscale basis functions:
\begin{description}
\item[General:] The type of partition one would use to upscale the model and/or form a coarse-scale discretization; that is, a partition with small variations in block sizes. Such partitions can be rectilinear/structured, or generated by a graph partitioning algorithm \citep{Karypis1998}, possibly with transmissibilities as the connection strength to minimize the permeability variation inside each block. In any configuration, our multiscale solvers will use at least one such partition with a prolongation operator generated from MsRSB basis functions.
\item[Static:] Partitions constructed to target specific features in the geocellular model. These partitions could be constructed by increasing the coarse-grid resolution near features of interest such as fractures and well paths, and/or by ensuring that block interfaces follow geological layers, fault surfaces, boundaries between different rock types, flow units, and depositional environments, etc. One effective way to generate such partitions is to agglomerate cells into blocks according to user-defined cell/face indicators and partitioning rules \citep{Hauge2012, Hauge2010, Lie2017general, lie2017mrstbook}. Several partition examples are shown in \cite{Lie2017}. We use MsRSB basis functions to define the resulting operators.
\item[Dynamic:] Basis functions designed to target dynamic couplings between pressure and saturations/mass fractions, i.e., the non-elliptic character of the pressure update. Such couplings are usually manifested as sharp transitions in the pressure updates across fluid interfaces. These partitions can be constructed using information (e.g., pressure update, velocity, saturations/mass fractions) from earlier nonlinear iterations. Herein, we either construct partitions by tracking displacement fronts, or by assigning cells into bins according to the magnitude of their pressure update from a previous nonlinear iteration. We use constant basis functions to construct the corresponding prolongation operators. 
\end{description}

We typically compute general and static partitions and corresponding basis functions during a preprocessing step. Basis functions can also be updated locally during the simulation to account for changing reservoir/fluid properties and driving forces by iterating a few times extra on the localized (elliptic) residual equations that define the MsRSB basis functions. The dynamic partitions must be recomputed during the simulation to account for the movement of fluids. This operation is not very computationally expensive, since it essentially consists of comparing two or more floating point numbers per cell and assigning each cell to a coarse block (i.e., setting an integer number in a partition vector). The optimal frequency for recomputing these partitions is obviously a trade-off between accurately capturing coupling between pressure and saturations/mass fractions and minimizing computational effort. Herein, we always use information (e.g., pressure update/phase saturations/mass fractions) from the most recent nonlinear iteration to compute dynamic basis functions before the next nonlinear iteration. We suggest the following update approach during a single time step:
\begin{enumerate}
    \item \emph{Before} the first nonlinear iteration, compute new dynamic basis functions based on information from the last nonlinear iteration of the previous time step.
    \item We also update the dynamic basis functions before the second nonlinear iteration of the time step.
    \item Before subsequent nonlinear iterations, we may also update the dynamic basis functions at a given update frequency.
\end{enumerate}
This update approach is motivated by the fact that during the solution to a time step, large regions of nonzero Newton updates are usually resolved within the first few iterations, and subsequent Newton updates are typically restricted to a small portion of the grid cells within the same regions \cite{Sheth2017}. Figure~\ref{fig:multiscale-bases} shows one basis function for each of the three types of partitions.

\subsection{Description of the solution procedure}

\begin{figure*}[t]
  \includegraphics[width = \textwidth]{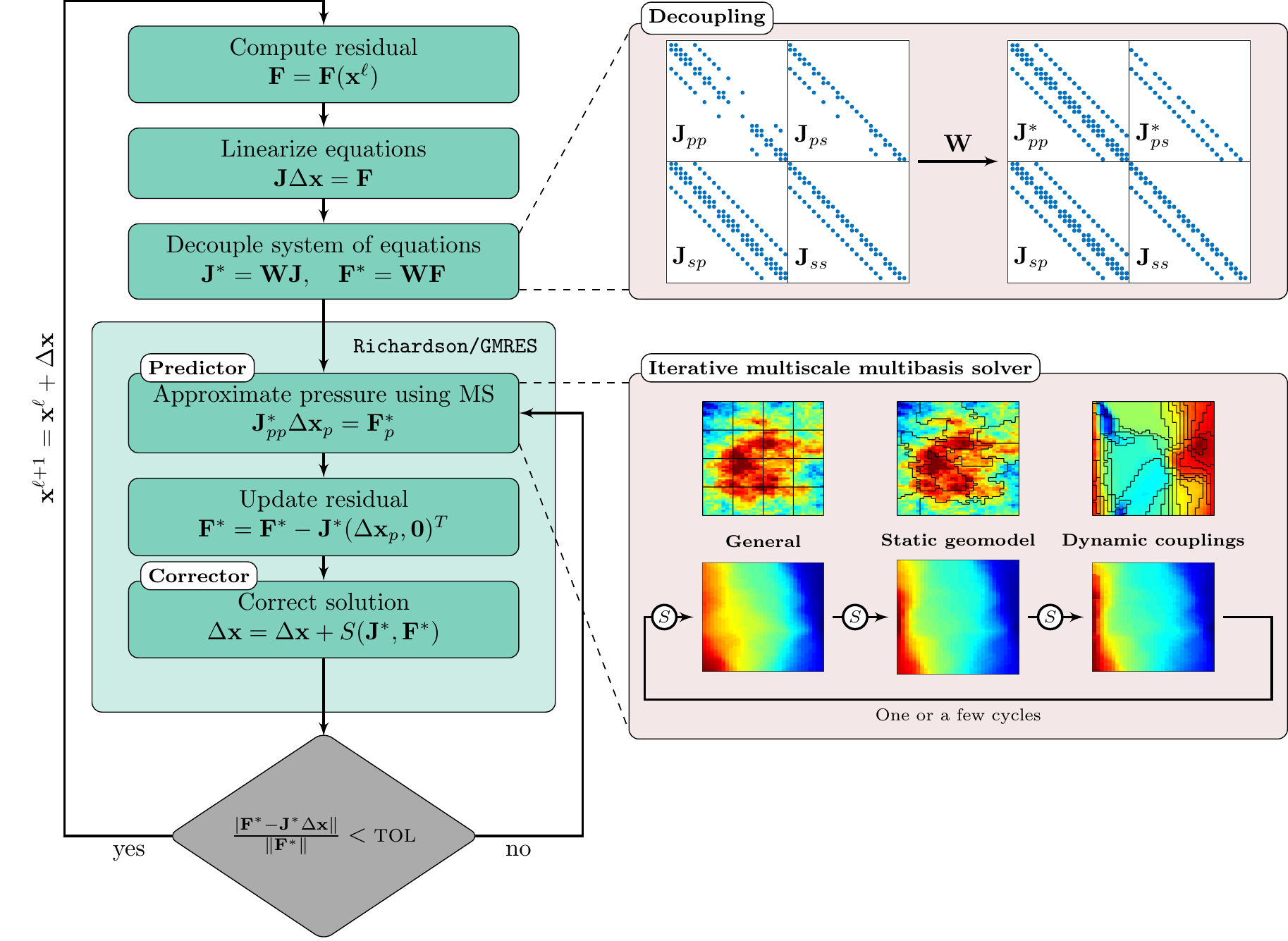}
  \caption{Flow chart describing the entire solution procedure to advance the solution one time step. After linearization, the linear system is decoupled. This gives a pressure block $\mat J_{pp}^*$ with a typical elliptic structure, and a coupling block $\mat J_{ps}^*$ with a typical hyperbolic structure, as seen in the \emph{Decoupling} box. The predictor step consists of one or a few cycles of the iterative multiscale multibasis method, as described in Equation~\eqref{eq:multiscale-multibasis}. To reduce high-frequency error modes, we apply a smoother between each multiscale solution, indicated by an $S$ in the \emph{Iterative multiscale multibasis solver} box.}
  \label{fig:solution-procedure}
\end{figure*}

One can combine the multiscale method with various iterative methods to form the pressure step of CPR. Herein, we either use a standard Richardson iteration, or use the multiscale solver as a two-step preconditioner in a Generalized Minimum Residual solver (GMRES) \citep{Saad1986}. The main steps in the solution procedure involved in one nonlinear iteration are thus: \emph{(i)} linearization to obtain the linear system, \emph{(ii)} decoupling to precondition the system for CPR, \emph{(iii)} the predictor step in which we approximate the pressure using an iterative multiscale multibasis solver, and \emph{(iv)} a correction step in which we smooth the solution to the full system. Figure~\ref{fig:solution-procedure} depicts the entire solution procedure.

\section{Numerical examples}
\label{sec:numerical-examples}
To test the effectiveness of the multiscale--multibasis method introduced above, we have implemented it in the open-source MATLAB Reservoir Simulation Toolbox (MRST) \citep{lie2017mrstbook, Krogstad2015}. Our conceptual implementation utilizes MATLAB quite efficiently, but is not yet fully optimized and has certain overhead one would not expect to see in a fully optimized code written in a compiled language. We therefore assess the computational efficiency in two steps: First, we present a number of test cases to assess how effective the multiscale, multibasis solver is as a CPR preconditioner. In each test case, we configure the approximate CPR solver to use various combinations of general, static, and dynamic partitions and then investigate which configuration requires fewest linear steps to reach the prescribed linear tolerance, averaged over all the nonlinear steps. Once we have investigated how the use of multiple multiscale bases affects the error reduction in the approximate CPR preconditioner, we analyze the theoretical number of floating-point operations required for the various stages of one linear iteration. Comparing the operational count for a multibasis preconditioner to a basic setup with a single MsRSB basis then gives us an indication of the potential for speedup.

\subsection{Example 1: Coupling strength}

\begin{figure*}[p!]
    \centering
    \includegraphics[width=0.9\textwidth]{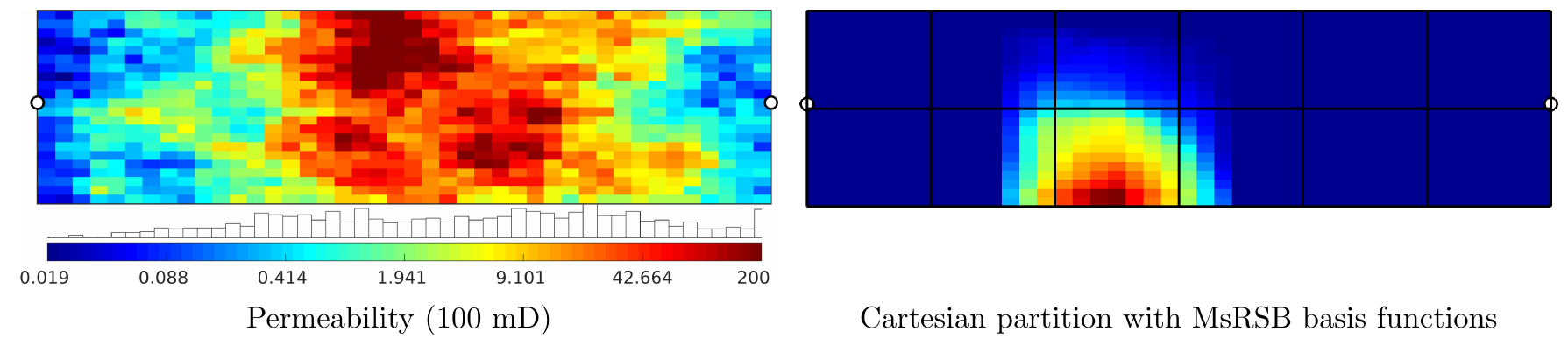}
    \caption{Permeability and general basis for Example~1.}
    \label{fig:coupling-setup}
    \vskip\floatsep
    \centering
    \hspace{-1em}\includegraphics[width=0.9\textwidth]{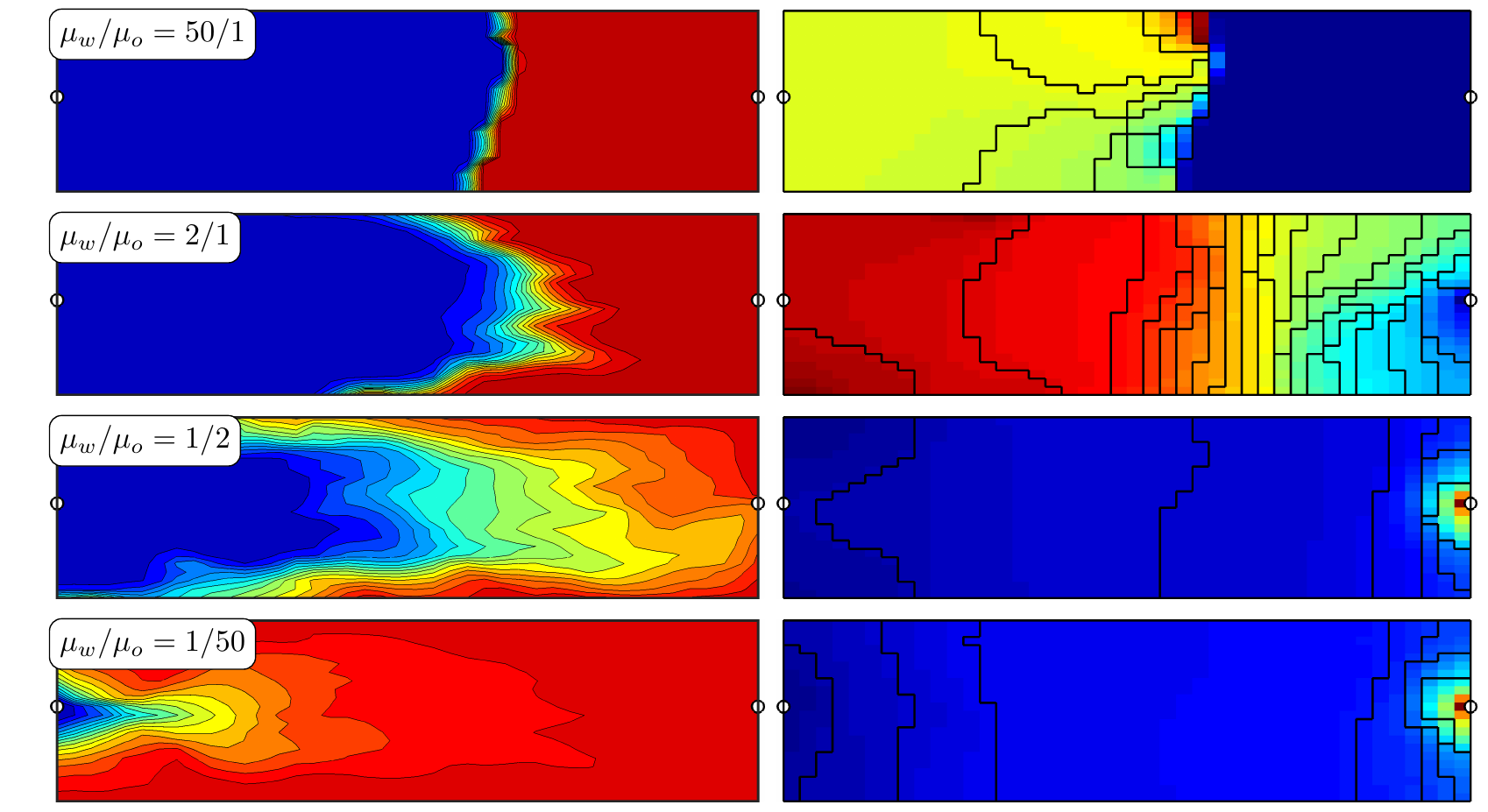}
    \caption{The left column shows saturation profiles after 900 days with four different viscosity ratios for Example~1. The right column shows the corresponding dynamic partitions (black lines), with colors indicating the magnitude of the pressure update.}
    \label{fig:coupling-visc-sat}
    \vskip\floatsep
    \def\fw{0.65\textwidth}
    \centering
    \includegraphics[width = \fw]{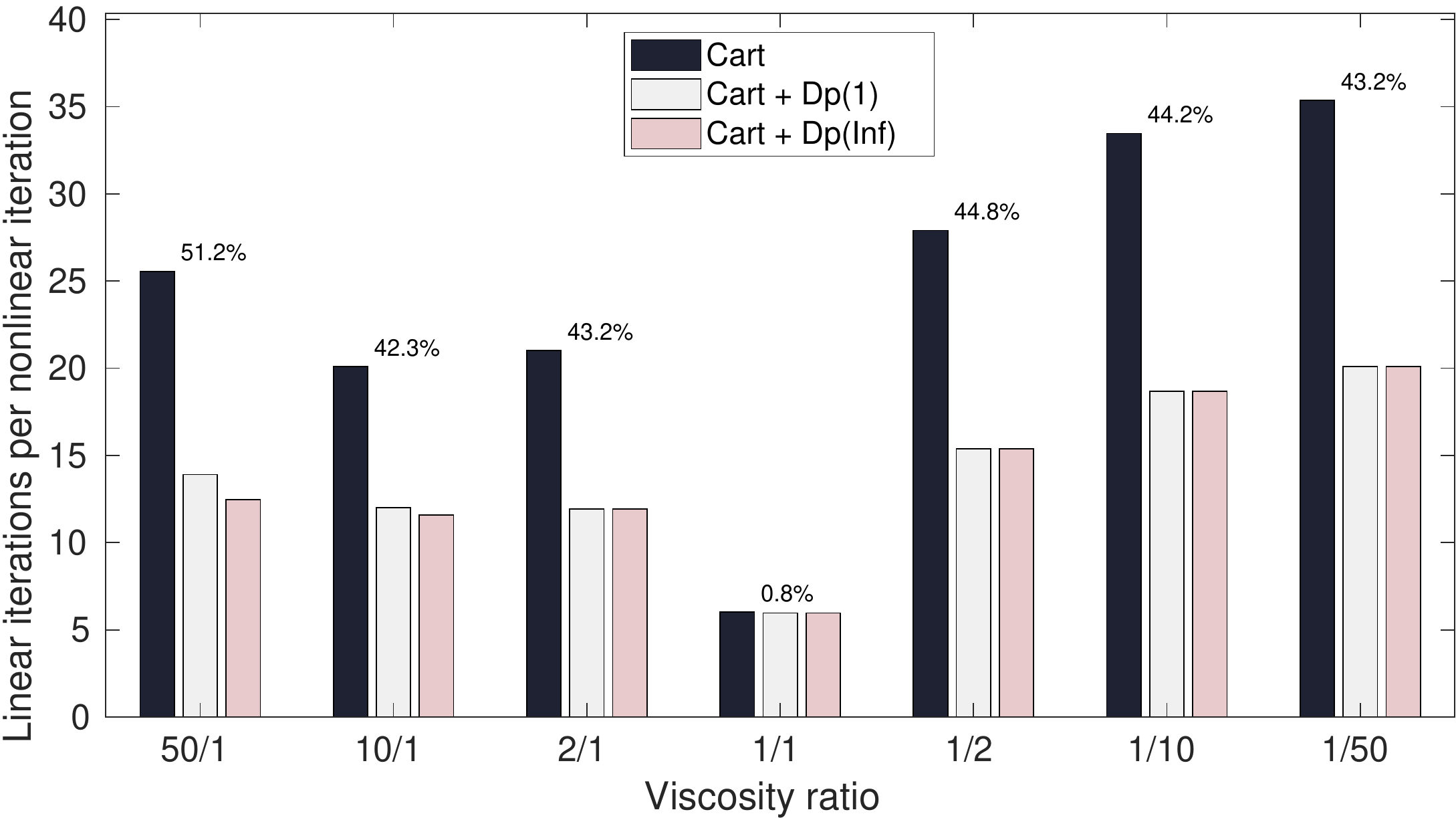}
    \caption{Average number of linear iterations per nonlinear iteration used by the three solvers for different viscosity ratios in Example~1. Dp(1) indicates recomputation of dynamic basis functions before each nonlinear iteration, whereas Dp(Inf) indicates recomputation before nonlinear iteration one and two. Percentages indicate the relative difference between the best and worst performing solver.}
    \label{fig:coupling-visc-its}
\end{figure*}

To study how the effect of using dynamic basis functions depends on coupling strengths, we consider a simple incompressible two-phase model with an aqueous and a liquid phase, each consisting of a single component. Referring back to \eqref{eq:incomp-elliptic}, we know that the coupling between pressure and saturation is determined by the viscosity ratio $\mu_\wat/\mu_\liq$, the density ratio $\rho_\wat/\rho_\liq$, and the magnitude of the capillary pressure $P_{c\liq\wat}$. The geological model consists of a $42\times 22$ subset from Layer 10 of Model 2 from the SPE10 benchmark \citep{spe10}, with the origin corresponding to cell $(10,110)$ in the full layer. The reservoir is initially filled with the liquid phase, and we inject one pore volume of the aqueous phase through an injection well at the left boundary over a period of four years. A producer operated a bottom-hole pressure of 50 bars is placed at the right boundary. A rectilinear $6 \times 2$ partition serves as the general basis (see Figure~\ref{fig:coupling-setup}), and we use the pressure update $\Delta p$ to generate dynamic basis functions. The incompressible flow equation \eqref{eq:incomp-elliptic} is already elliptic and could be solved directly without any CPR reduction, but herein we use the general procedure outlined above. In our experience, multiscale solvers based on Richardson iterations are more sensitive to the quality of the preconditioner than solvers based on GMRES. We thus use a Richardson-type solver to emphasize the effect of adding a dynamic basis.

\begin{figure*}[t!]
    \centering
    \includegraphics[width=0.9\textwidth]{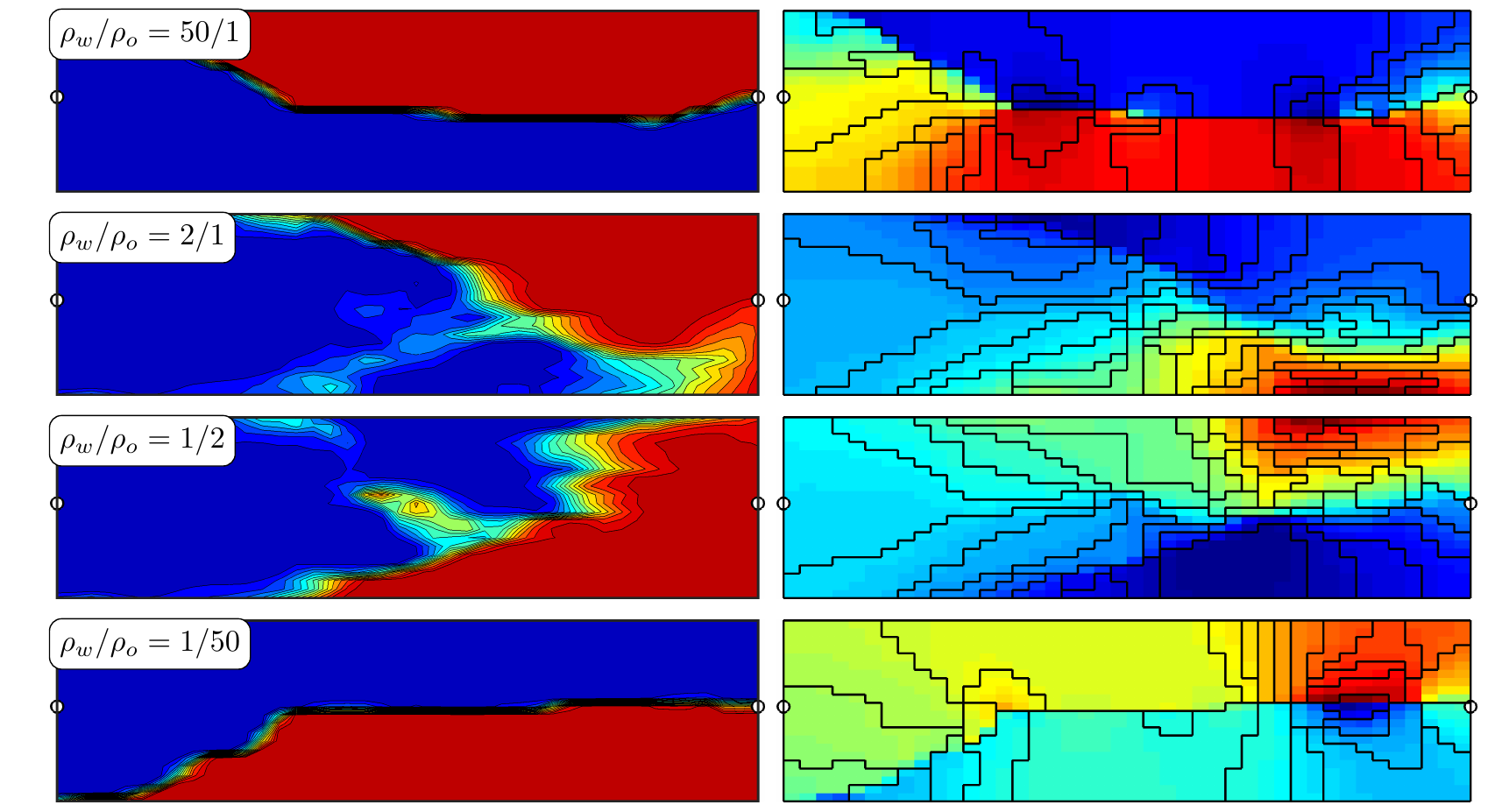}
    \caption{The left column shows saturation profiles after 900 days for four different density ratios in Example~1. The right column shows the corresponding dynamic partitions (black lines), with colors indicating the magnitude of the pressure update.}
    \label{fig:coupling-grav-sat}
    \vskip\floatsep
    \def\fw{0.65\textwidth}
    \centering
    \includegraphics[width = \fw]{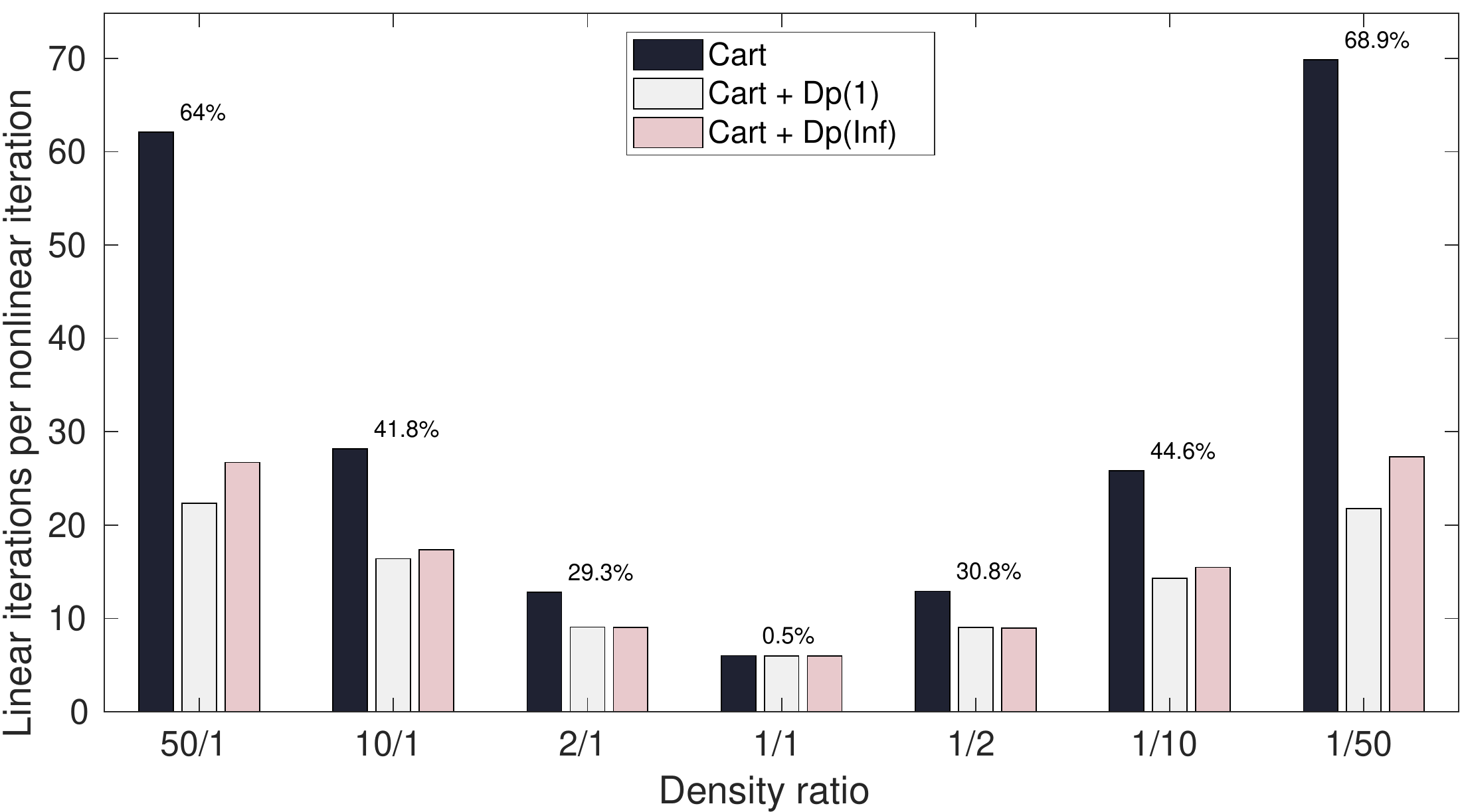}
    \caption{Average number of linear iterations per nonlinear iteration used by the three solvers for different density ratios in Example~1. Dp(1) indicates that we recompute dynamic basis functions before each nonlinear iteration, whereas Dp(Inf) indicates that we only recompute before nonlinear iteration one and two. Percentages indicate the relative difference between the solvers with best and worst performance.}
    \label{fig:coupling-grav-its}
\end{figure*}

We neglect gravity and capillary forces, assume linear relative permeabilities for both phases, and let the viscosity ratio $\mu_\wat/\mu_\liq$ vary from $50/1$ to $1/50$, normalized so that the less viscous phase has a viscosity of 1 cP. Figure~\ref{fig:coupling-visc-sat} shows saturation profiles after 900 days for selected viscosity ratios. We clearly see how the dynamic partitions shown in the right column of the figure adapt to the saturation front for the two cases where the injected fluid is more viscous. Figure~\ref{fig:coupling-visc-its} reports the average number of linear iterations used for the different solvers, i.e., the number of iterations performed internally in the Richardson/GMRES box in Figure~\ref{fig:solution-procedure} averaged over all the nonlinear iteration steps. We see that adding a dynamic basis is beneficial in all cases except with unit viscosity ratio. This is because pressure and saturation in this case are completely decoupled, and the negligible difference of 0.8 \% comes only from the smoothing iteration in the first step of \eqref{eq:multiscale-multibasis}. We have also compared the effect of recomputing the dynamic basis functions before each nonlinear iteration instead of only before the first and second iterations. The solver with update before each iteration is denoted Dp(1), whereas the solver with update before iteration one and two is denoted Dp(Inf). Interestingly, we observe no reduction in the number of iterations by recomputing before each time step, and even an increase for large viscosity ratios. A possible explanation to this is that the last Newton updates of each time step are very small in most of the domain, so that the construction of dynamic partitions is influenced by noise.

\begin{figure*}[t!]
    \centering
    \includegraphics[width=0.95\textwidth]{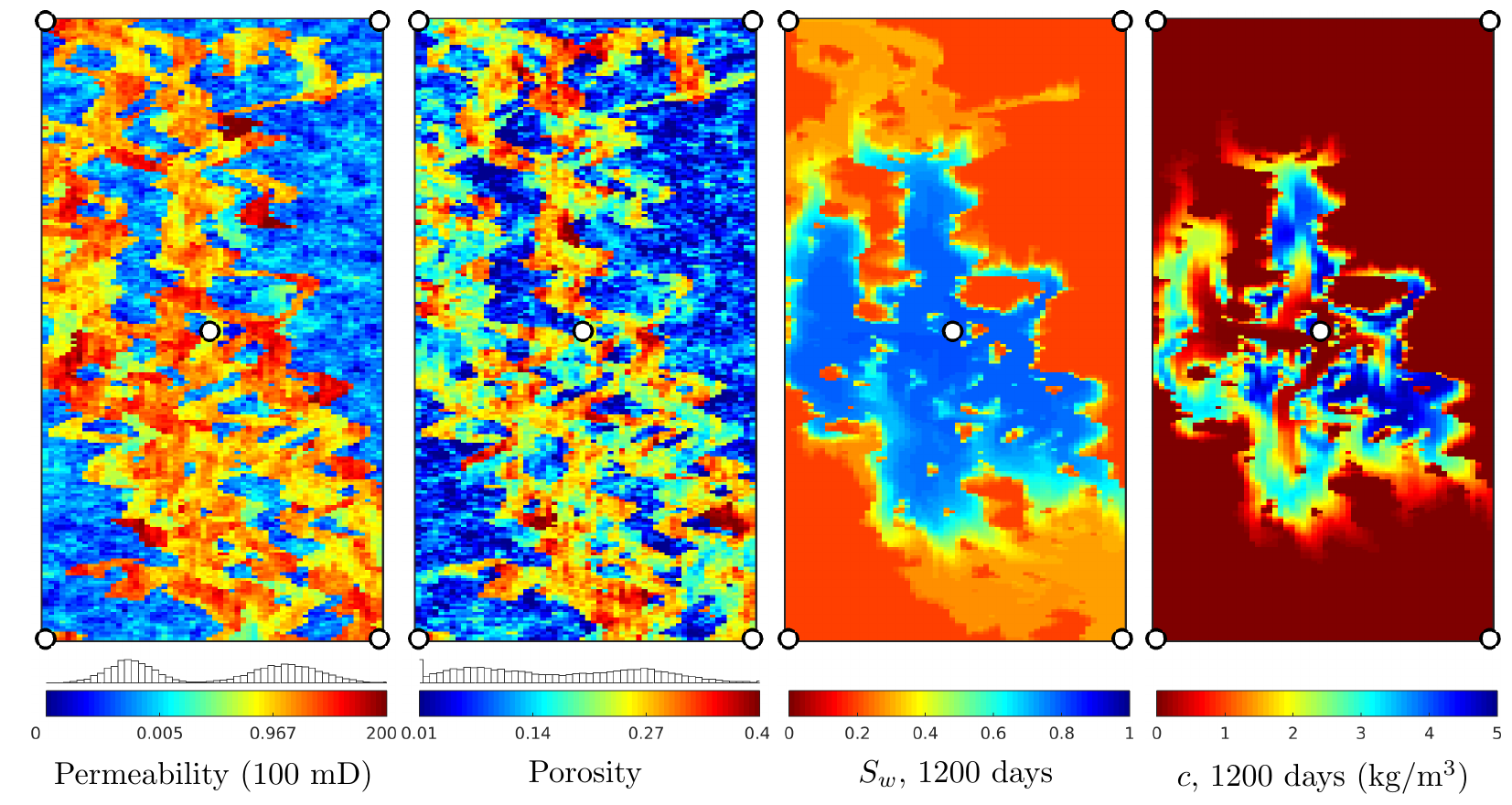}
  \caption{Petrophysical properties for the case in Example~2, along with water saturation and polymer concentration after 1200 days of injection.}
  \label{polymer-setup}
\end{figure*}

Next, we set equal viscosities and look at the effect of varying the density ratio in the presence of gravity; that is, we vary the ratio $\rho_\wat/\rho_\liq$ from $50/1$ to $1/50$, normalized so that the less dense fluid has a density of 10 kg/m$^3$. Gravity is aligned with the negative $y$-direction. Figure~\ref{fig:coupling-grav-sat} shows the saturations after 900 days for different density ratios. We clearly see how the dynamic partitions adapt to the interface between the segregated fluids, especially for the highest and lowest density ratios. Figure~\ref{fig:coupling-grav-its} reports the  average number of linear iterations per nonlinear iteration used by the two solvers. Again, we see that adding a dynamic basis improves the convergence rate of the linear solver except for the case with unit density ratio. Moreover, compared with the case with varying viscosity ratios, the relative reduction in iterations is bigger and depends more strongly on the density ratio. For the extreme density ratios of $\rho_\wat/\rho_\liq = 50/1$ and $1/50$, the gravity segregation will take place very rapidly, giving sharp fluid interfaces where the pressure update is nearly discontinuous. These interfaces are nicely captured by the dynamic basis functions. In this case, we observe a slight reduction in the number of linear iterations by recomputing the dynamic basis functions before each nonlinear iteration (Dp(1)).

\subsection{Example 2: Polymer injection}

We consider polymer injection in Layer 52 of SPE10 Model 2 with wells placed in an inverted quarter five-spot pattern. The center well injects water at a constant rate of 9.14 m$^3$/day over a period of 2000 days, and the four corner wells operate at fixed bottom-hole pressures of 275 bar. The injected water will rapidly fill the high-permeable fluvial channels. To deviate water into the low-permeable zones, we inject a polymer slug (5 kg/m$^ 3$) from 400 to 800 days. The polymer is assumed to mix completely with water \citep{ToddLongstaff1972}, increasing the viscosity of the injected phase from 0.3 cP to 60 cP. The polymer model also accounts for adsorption of polymer onto the porous rock. Shear effects and inaccessible pore space can also be accounted for in the model, but this is for simplicity not included in the example. See \citep{Bao2017} for further details. The polymer effect on viscosity introduces a strong coupling between pressure and water saturation/polymer concentration. Figure~\ref{polymer-setup} shows petrophysical properties, well positions, water saturation and polymer concentration after 1200 days. Unlike in the previous example, the system has compressibility and we thus use a CPR solver.

\begin{figure*}[t!]
    \centering
    \includegraphics[width=\textwidth]{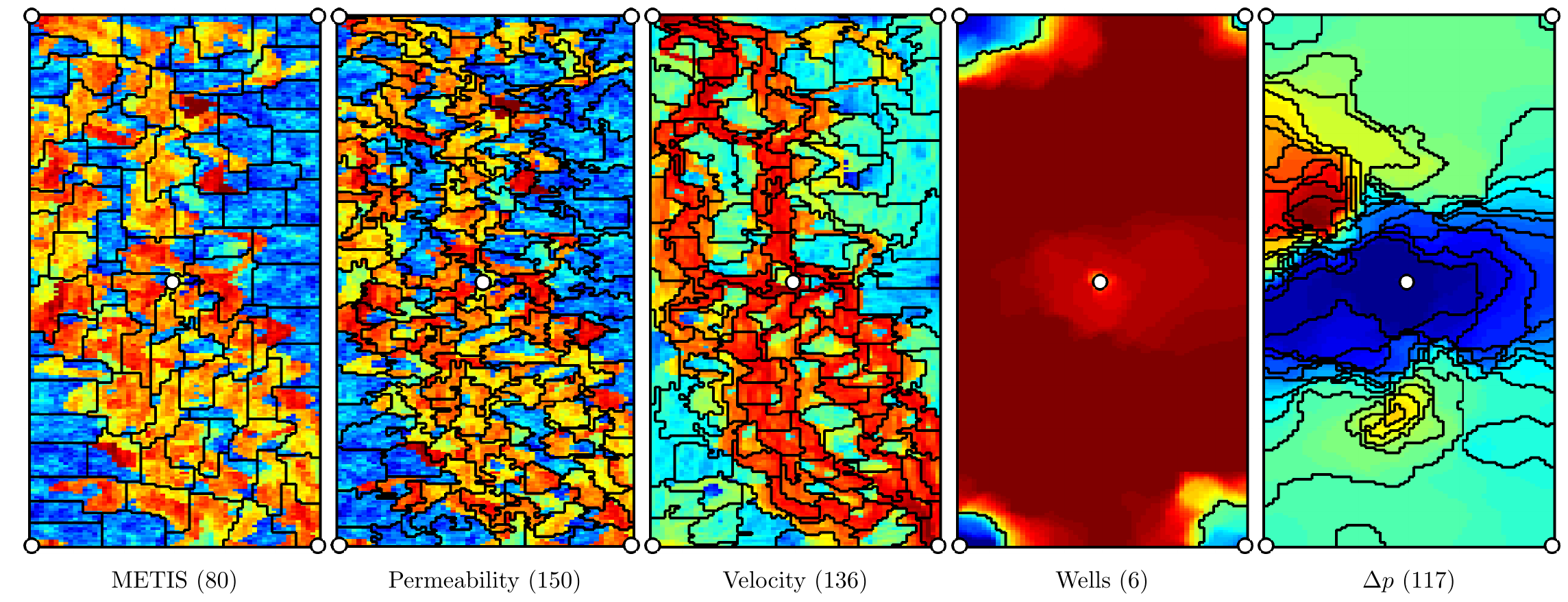}
    \caption{Partitions for the different basis functions used in the polymer flooding in Example~2. The number of coarse cells for each partition are indicated in parentheses. The colors show the indicator used to generate the different partitions, with the exception of the well partition, where the actual basis functions are shown.}
    \label{fig:polymer-partitions}
    \vskip\floatsep
    \centering
    \includegraphics[width = 0.9\textwidth]{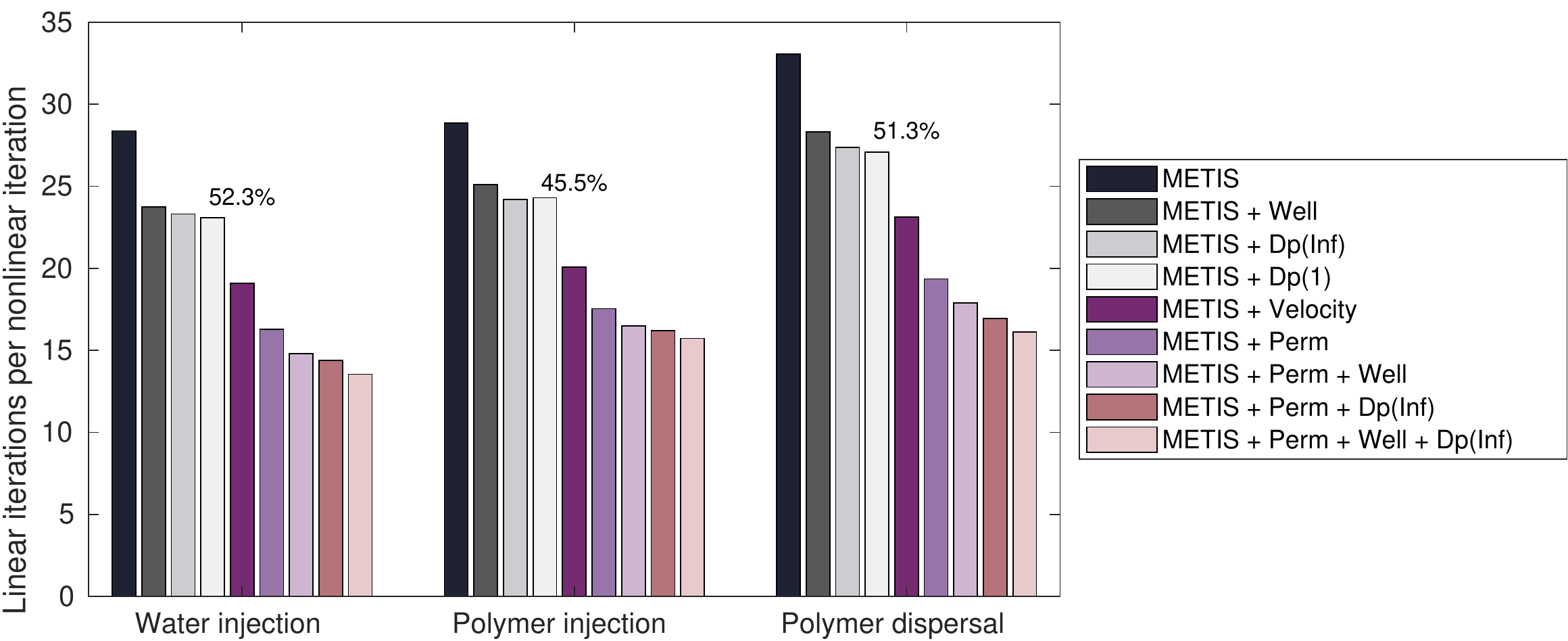}
    \caption{Average number of linear iterations per nonlinear iteration with different combinations of multiscale operators for the CPR preconditioner in Example~2. Dp(1) indicates recomputation of dynamic basis functions before each nonlinear iteration, whereas Dp(Inf) indicates recomputation before nonlinear iteration one and two. Percentages indicate the relative difference between the best and worst performing solver.}
\label{fig:polymer-bar}
\end{figure*}

Figure~\ref{fig:polymer-partitions} shows the different types of partitions we use to generate basis functions for the CPR preconditioner. The general basis consists of MsRSB basis functions defined on a partition constructed by METIS \citep{Karypis1998} with one-sided transmissibilities as edge weights. We also consider two partitions that honor the high permeability contrast and channeled structure of the reservoir, generated by an ad-hoc algorithm that agglomerates cells into coarse blocks based on a flow indicator \citep{Hauge2010}. The first partition uses permeability as flow indicator, whereas the second uses an incompressible velocity field computed a priori for the same grid and well setup. For both we use MsRSB basis functions to construct the corresponding prolongation operator. In addition, we use well bases and a set of dynamic basis functions with the pressure update $\Delta p$ from the previous Newton iteration as indicator to agglomerate the grid cells.

\begin{figure*}[t!]
    \centering
    \includegraphics[width=0.95\textwidth]{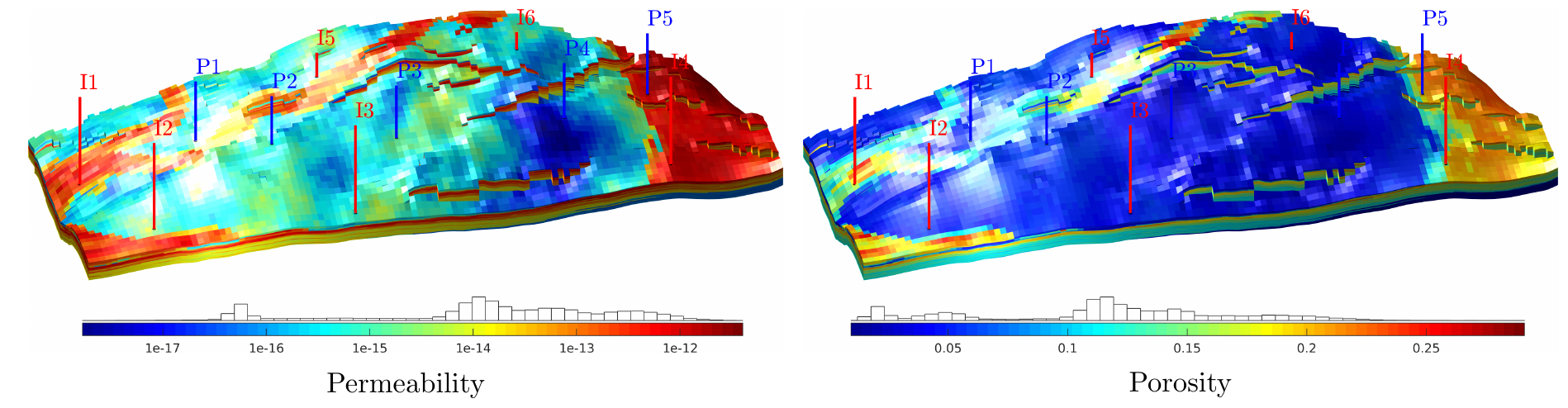}
    \caption{Petrophysical properties and well configuration for the SAIGUP field in Example 3.}
    \label{fig:saigup-setup}
    \vskip\floatsep
    \includegraphics[width=0.95\textwidth]{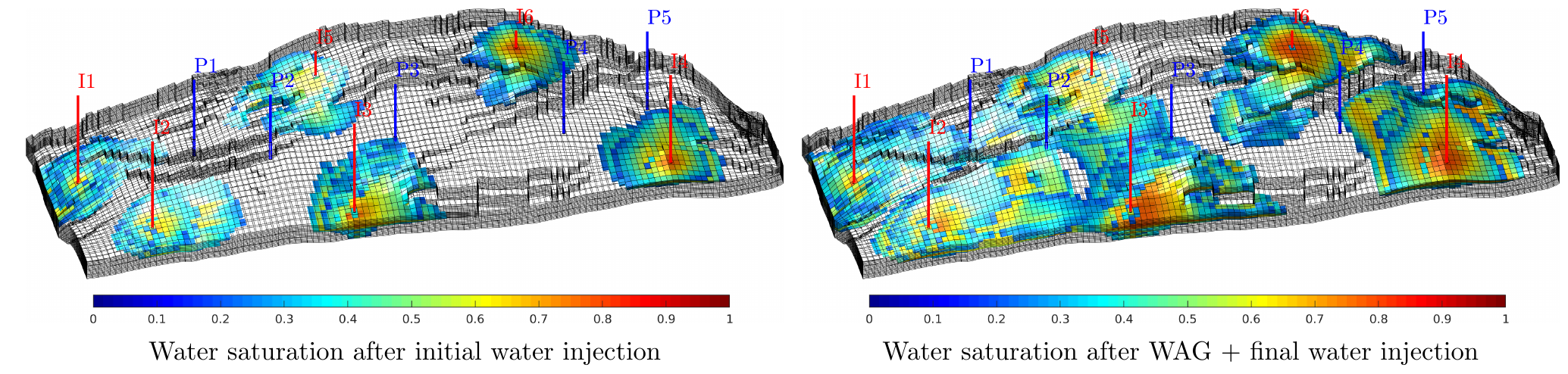}
    \caption{Water saturation before and after WAG + water injection in Example 3.}
    \label{fig:saigup-sat}
\end{figure*}

Figure~\ref{fig:polymer-bar} reports the average number of linear iterations per nonlinear iteration for the different solvers during three separate phases of the recovery process: initial water injection, injection of the polymer slug, and the final water-injection period that disperses the injected polymer slug. The most significant reduction in iterations comes from adding static partitions that adapt to the high-flow channels, identified either by permeability or single-phase velocities. Using well basis functions typically gains 3--4 iterations. This comes at a very low computational cost, since the basis only consists of six coarse blocks (one for each well plus one to ensure partition of unity). The dynamic partitions have a similar reduction effect, but are generally more costly than the well bases. The dynamic bases reduce the number of iterations a bit more during the polymer dispersal phase, when the effective viscosity ratio between the displaced and displacing fluid is at its highest. Similar to the previous example, we observe only a slight reduction in the number of linear iterations when recomputing the dynamic basis functions before each nonlinear iteration (Dp(1)).

\subsection{Example~3: Field model (SAIGUP)}
The SAIGUP study \citep{Manzocchi2008} developed several realistic models of shallow-marine oil reservoirs based on data from real shoreface deposition sites, here represented by a $40\times 120\times 20$ corner-point grid that spans a lateral area of approximately $9\times 3$ km$^2$. The model has 78,720 active cells and describes several major faults. Permeability and porosity follow multimodal distributions, and the model identifies six different rock types with clearly distinct petrophysical characteristics. The reservoir also contains large amounts of mud-drapes that reduce the vertical permeability. Figure~\ref{fig:saigup-setup} shows porosity and  lateral permeability and the location of the injection and production wells.

We consider a water-alternating-gas (WAG) injection with six injectors operating at a constant rate and five producers operating at constant bottom-hole pressure of 50 bar. WAG is a popular enhanced oil recovery strategy in which small volumes of water and gas are injected in a cycle. The injected gas will dissolve into the reservoir oil, which improves sweep efficiency and releases much of the trapped residual oil. Water is injected between the gas volumes to uphold the mobility of the displacement front. This technique is typically employed after an initial water flood. In our case, we start by injecting 0.8 pore volumes of water over a period of two years. Over the next two years we inject 0.8 pore volumes of alternating water and gas, followed by another 0.8 pore volumes of water. 

\begin{figure*}[t!]
    \centering
    \includegraphics[width=0.95\textwidth]{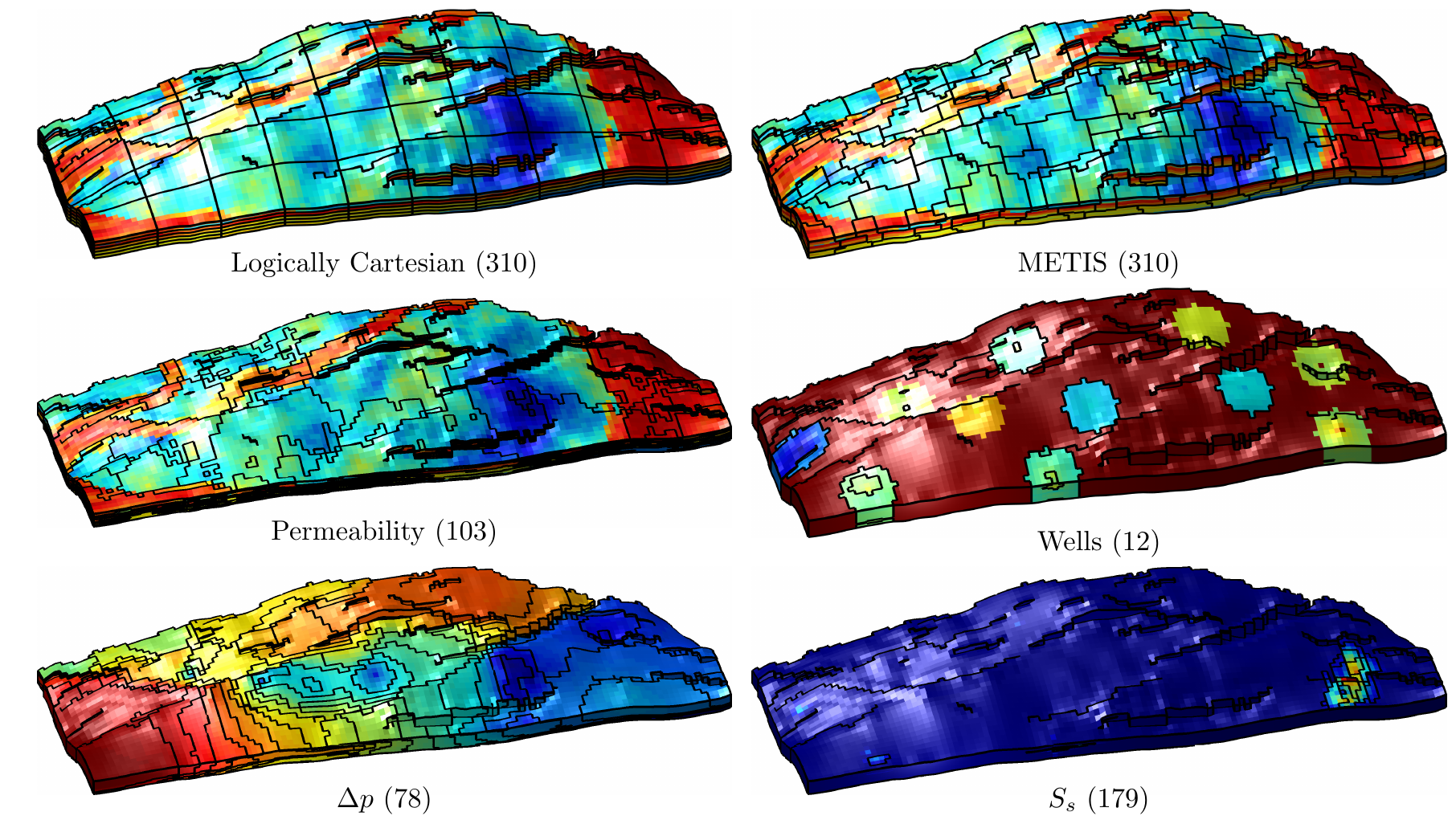}
    \caption{The six different partitions used for the WAG scenario in Example 3. The number of coarse blocks for each partition is indicated in parentheses. Colors show the indicator used to generate each partition, except for the well partition, where we show the actual basis functions.}
    \label{fig:saigup-grids}
    \vskip\floatsep
    \includegraphics[width = 0.9\textwidth]{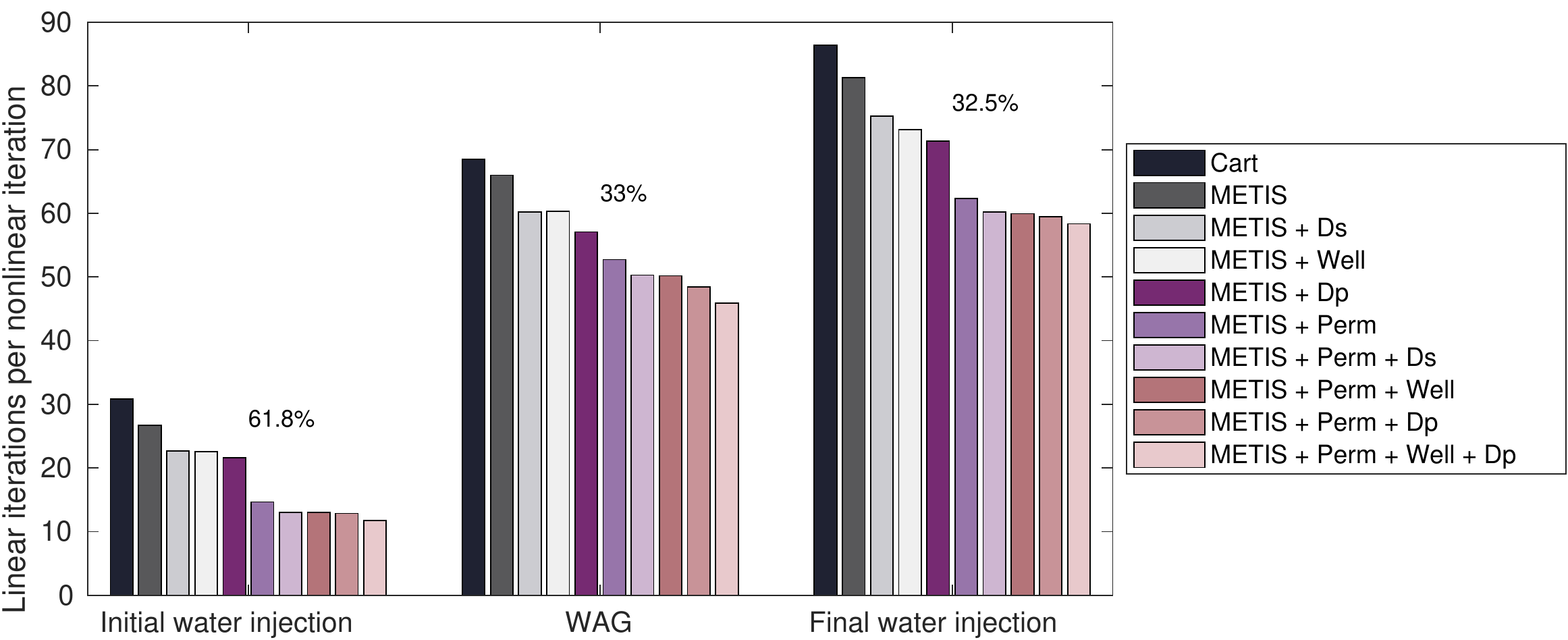}
    \caption{Number of linear iterations per nonlinear iteration used by the different CPR solvers during each of the three production stages in Example~3. Percentages indicate the relative difference between the best and worst performing solver.}
    \label{fig:saigup-bar}
\end{figure*}

Oil and gas are slightly compressible, with quadratic Brooks-Corey relative permeabilities. Densities  and viscosities are 1000, 800, and 100 kg/m$^3$ and 1, 3, and 0.4 cP at reference pressure for the water, oil, and gas phases, respectively. The injected solvent gas has a density of 100 kg/m$^3$ and a viscosity of 0.5 cP, and mixes with the reservoir oil and gas according to a modified version of the Todd--Longstaff mixing rule \citep{ToddLongstaff1972}. This model is the same as implemented in a commercial simulator \citep{eclipse}. With only reservoir oil and reservoir gas, we have a traditional black-oil behavior, whereas the reservoir oil is assumed to be fully miscible with the solvent gas when no reservoir gas is present. In the intermediate region, we interpolate between the two extremes. The oil/solvent mixture has much lower density and viscosity than pure oil. In addition, solvent gas is assumed to lower the residual oil saturation from $S_{or,i}=0.2$ in the immiscible case to $S_{or,m}=0.1$ in the fully miscible case. The strong coupling between solvent gas saturation and reservoir fluid mobility means that the solution is strongly affected by the solvent saturation. As a consequence, the two-stage CPR preconditioner may not be as effective as for a standard waterflooding scenario. Figure~\ref{fig:saigup-sat} reports water saturations after initial water injection and at the end of the simulation. 

Figure~\ref{fig:saigup-grids} shows the multiscale partitions used for the CPR preconditioner. These include a logically Cartesian partition and a partition generated using one-sided transmissibilities to measure coupling strengths between cells in the METIS graph-partitioning algorithm. These serve as general partitions. We include a static partition to honor the large variations in permeability and another static well partition that accounts for changes in the well control. Finally, we use two dynamic partitions agglomerated with the pressure update and the solvent gas saturation from the previous nonlinear iteration as indicators. Based on observations from the preceding examples, we only recompute the dynamic basis functions after nonlinear iteration one and two in this example. Figure~\ref{fig:saigup-bar} reports average number of linear iterations per nonlinear iteration for different combinations of multiscale operators. All solvers consume significantly more iterations during the WAG and the final water injection period. The presence of solvent gas gives denser systems, amplified by the fact that the coupling between pressure and solvent saturation increases as the solvent gas sweeps the reservoir and mixes with the resident fluids. 

Using METIS to group cells with similar permeabilities gives slightly faster convergence than a Cartesian partition in index space. Adding extra static and/or dynamic partitions reduces the iterations further, and in this particular example, adapting the static partition to permeability is the most effective. Trends are similar during all three stages of the production period, but the relative differences between the solvers are smaller during WAG and final water injection than during the initial water injection.

\begin{figure*}[t!]
    \centering
    \hspace{-2em}\includegraphics[width=0.95\textwidth]{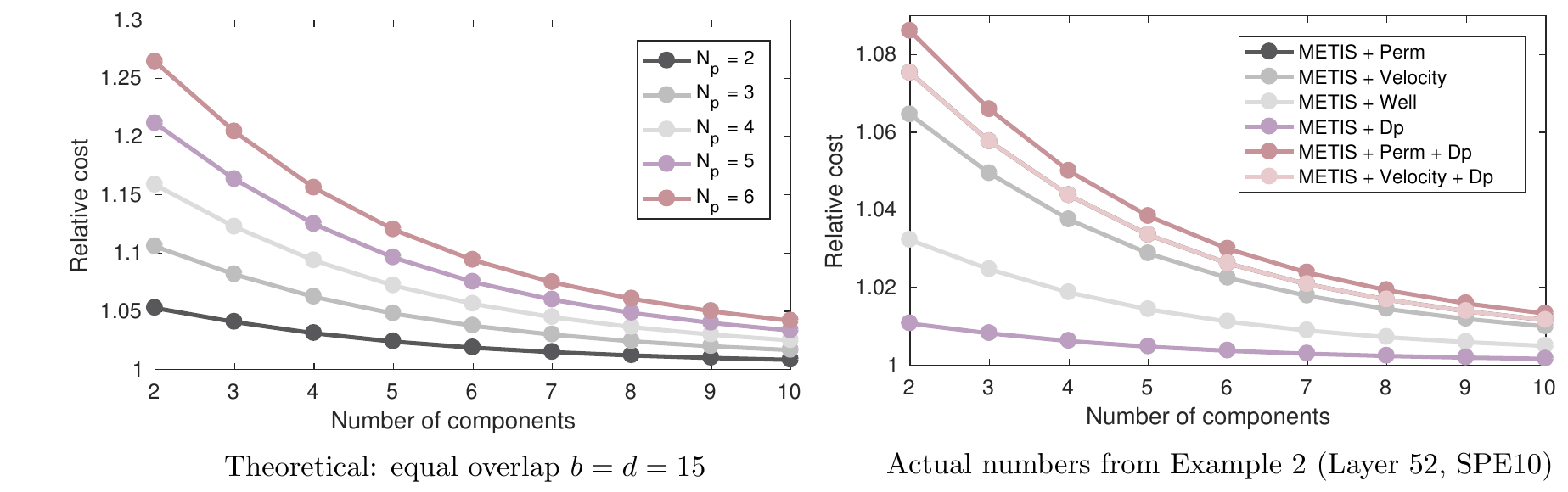}
    \caption{The relative added computational cost for a full linear step using multiple instead of a single multiscale basis in the approximate CPR preconditioner reported as function of the number of components, i.e., primary unknowns per cell.}
    \label{fig:complexity}
\end{figure*}

\subsection{Computational efficiency}

To assess the computational costs of the methods discussed above, we compare the theoretical number of operations required for one full linear iteration. From Equation~\eqref{eq:multiscale-multibasis}, we see that this involves the following steps:
\begin{description}
    \item[Smooth:] We use incomplete LU-factorization with zero fill-in (ILU(0)) as our smoother. With one pressure unknown for each of the $\ncell$ cells, the smoothing step consists of solving two triangular systems and consumes $\mathcal{O}(2\ncell\nsparse)$ operations, where $\nsparse\ll \ncell$ is an upper bound on the number of nonzero elements in each row of $\mat A$. Computing the linearized residual, $\mat q - \mat A\vec x$, amounts to $\ncell \nsparse+N$ operations for the matrix-vector product and the subtraction. Updating the state vector $\vec x$ adds another $\ncell$ operations.
    \item[Restrict:] This amounts to computing $\vec d_\ell = \mat R^\ell(\mat q - \mat A \vec x)$. The associated cost is $\mathcal{O}(\ncell\nsparse+N)$ operations to compute the linear residual and $\mathcal{O}(\noverlap_\ell\ncell)$ operations for multiplication by $\mat R^\ell$. Here, $\noverlap_\ell=1$ for a finite-volume restriction operator, whereas $1<\noverlap_\ell\ll\nblock_\ell$ is the maximum number of basis functions with support in a single cell for a Galerkin restriction.
    \item[Solve for $\vec p_c$:] Inverting the coarse matrix $(\mat A_c^\ell)^{-1}\vec d$ typically takes $\mathcal{O}(\nblock_\ell^p)$ operations with $p\approx 2$. The worst case is $p=3$ for a dense matrix inverted by direct Gaussian elimination.
    \item[Prolongate and update:] Prolongating the coarse pressure back to the fine grid, $\mat P^\ell \vec p_c$, consumes $\mathcal{O}(\noverlap_\ell \ncell)$ operations and updating the state $\vec x$ adds another $\ncell$ operations.
\end{description}
To summarize, we assume for simplicity that each of the $\npart$ coarse partitions consists of $\nblock$ blocks and that all basis functions have the same maximum overlap $\noverlap$.  Then, we have that one cycle involves the following number of operations
\begin{equation*}
    \begin{aligned}
        & \npart\mathcal{O}\bigl((\ncell\nsparse + \ncell + 2\ncell\nsparse + \ncell) + (\ncell\nsparse + \ncell + \noverlap\ncell) + \nblock^p \\
        & \qquad \qquad \qquad + (\noverlap\ncell + \ncell) \bigr) \approx \mathcal{O}\bigl((4\nsparse+2\noverlap+4)\npart\ncell\bigr)
    \end{aligned}
\end{equation*}
if we assume that the cost of inverting the multiscale system is negligible or of order $\mathcal{O}(\ncell)$.

Updating the residual after the predictor step involves multiplication of the approximated pressure $\Delta \vec x_p$ by the first $\ncell$ columns and all $\ncell\ncomp$ rows of $\mat J^*$ (in which each row will have at most $\nsparse$ nonzero elements) and subtracting the result from the residual. All in all, this involves $\mathcal{O}(\ncomp\ncell\nsparse$ + $\ncomp\ncell)$ operations. Finally, the corrector step consists of an ILU(0) smoothing step applied to the full system with $\ncell\ncomp$ unknowns for a model with $\ncomp$ components. This system will also have approximately $\ncomp$ as many nonzero elements on each row. Updating the result gives $\ncomp\ncell$ more operations. Altogether, this amounts to a total of $\mathcal{O}\left(2\nsparse\ncomp^2\ncell\right) + \ncomp\ncell$ operations.

To sum up, given that we perform only one cycle in the iterative multiscale multibasis solver, the cost of one linear iteration using $\npart$ multiscale partitions is approximately
\begin{equation}
\begin{split}
    c(\npart) =  \ncell\bigl[& \npart(4\nsparse + 2\noverlap + 4) \\ & + \ncomp(\nsparse + 1) + \ncomp(2\nsparse\ncomp + 1)\bigr].
    \end{split}
\end{equation}
Figure~\ref{fig:complexity} exemplifies the complexity analysis by reporting the ratio between the cost $c(\npart)$ of using $\npart$ different multiscale partitions and the cost $c(1)$ of using only one multiscale partition for two different cases. In the first case, we assume that all $\npart$ coarse partitions have the same maximum overlap $\noverlap$, and that $\nsparse = 15$, which is not uncommon for a realistic reservoir model with many faults and other types of non-neighboring connections. The second case uses numbers from Example~2 discussed above. We see that the worst-case scenario consists of problems with few components (primary unknowns per cell) and highly unstructured grids with many cell neighbors. As expected, we also see that the relative cost of using a multibasis CPR preconditioner diminishes with increasing number of components. 

\section{Concluding remarks}
We have reported a preliminary study of how multiple operators that adapt to static and/or dynamic features can improve the efficiency of multiscale-CPR preconditioning in fully implicit simulators. A number of numerical experiments show that one can reduce the total number of linear iterations required to solve the full discrete system by applying extra multiscale operators to the reduced pressure system. 

Adapted multiscale operators typically have few degrees of freedom relative to operators corresponding to general partitions. Applying one more iteration with such an operator to the reduced pressure system is therefore inexpensive compared with the cost of a linear iteration for the full system. From our simple analysis of the operational count, we conclude that the additional cost of using multiple bases for a typical scenario should be at most 10--15\%. 

Basis functions can be adapted and combined in many different ways, and our results are somewhat inconclusive with respect to which specific combination is the most effective.  In general, we expect that the best choice of multiscale operator(s) will vary from case to case and depend on the reservoir geology (level and type of heterogeneity), the drive mechanisms, and the degree of coupling between pressure and the other variables. Our experiments nevertheless give strong indication that it is beneficial to add static partitions honoring large permeability contrasts and/or add basis functions that better resolve near-well regions. Likewise, it also seems beneficial to add dynamic partitions adapting to pressure updates whenever these are located along propagating fluid fronts and/or fluid phase interfaces. Our overall conclusion is that using a multibasis CPR preconditioner may give a 10--60\% speedup compared with using a multiscale-CPR preconditioner with only a single set of basis functions.

Finally, we note that the idea of using dynamic basis functions may be beneficial in a predictor-corrector solution framework where the system is first solved sequentially in a predictor step, and then corrected by solving the fully implicit system using the sequential solution as an initial guess, see e.g., \cite{MM18:ecmor}. In this case, the sequential pressure solution will typically be close to the fully implicit solution, and dynamic basis functions based on this pressure may be very well suited for a multiscale pressure solver.

\section{Acknowledgements}
The work of Klemetsdal and Lie is funded by the Research Council of Norway through grant no.~244361. M{\o}yner is funded by VISTA, a basic research program funded by Equinor and conducted in close collaboration with the Norwegian Academy of Science and Letters.

\begin{small}
  \bibliography{paper.bib}
\end{small}

\end{document}